\documentclass{amsart}
\def\Defined{}
\ifx\UseRussian\Defined
	\usepackage[T2A,T2B]{fontenc}
	\usepackage[cp1251]{inputenc}
	\usepackage[english,russian]{babel}
\fi
\usepackage{footmisc}
\usepackage[all]{xy}
\usepackage{color}
\definecolor{UrlColor}{rgb}{.9,0,.3}
\definecolor{SymbColor}{rgb}{.4,0,.9}
\definecolor{IndexColor}{rgb}{1,.3,.6}
\definecolor{eml1}{rgb}{.8,.1,.1}
\definecolor{eml2}{rgb}{.1,.6,.6}

\usepackage{xr-hyper}
\usepackage[unicode]{hyperref}
\hypersetup{pdfdisplaydoctitle=true}
\hypersetup{colorlinks}
\hypersetup{citecolor=UrlColor}
\hypersetup{urlcolor=UrlColor}
\hypersetup{pdffitwindow=true}
\hypersetup{pdfnewwindow=true}
\hypersetup{pdfstartview={FitH}}

\def\hyph{\penalty0\hskip0pt\relax-\penalty0\hskip0pt\relax}
\def\Hyph{-\penalty0\hskip0pt\relax}

\newcommand{\Basis}[1]{\overline{\overline{#1}}{}}
\newcommand{\Vector}[1]{\overline{#1}{}}
\newcommand{\gi}[1]{\boldsymbol{\textcolor{IndexColor}{#1}}}
\makeatletter
\newcommand{\NameDef}[1]{%
	\expandafter\gdef\csname #1\endcsname%
}%
\newcommand{\ShowSymbol}[1]{%
	\@nameuse{ViewSymbol#1}%
}%
\newcommand{\symb}[3]{%
	\@ifundefined{ViewSymbol#3}{%
		\NameDef{ViewSymbol#3}{\textcolor{SymbColor}{#1}}%
		\NameDef{RefSymbol#3}{\pageref{symbol: #3}}%
		\@namedef{LabelSymbol#3}{\label{symbol: #3}}%
	}{%
		\NameDef{RefSymbol#3}{}%
		\@namedef{LabelSymbol#3}{}%
	}%
	\ifcase#2
	\or
		$\@nameuse{ViewSymbol#3}$%
	\or
		\[\@nameuse{ViewSymbol#3}\]%
	\else%
	\fi%
	\@nameuse{LabelSymbol#3}%
}%
\makeatother

\newcommand{\subs}{${}_*$\Hyph}
\newcommand{\sups}{${}^*$\Hyph}

\newcommand{\CRstar}{{}_*{}^*}
\newcommand{\RCstar}{{}^*{}_*}

\newcommand{\RC}{$\RCstar$\Hyph}
\newcommand{\CR}{$\CRstar$\Hyph}
\newcommand{\drc}{$D\RCstar$\Hyph}
\newcommand{\dcr}{$D\CRstar$\hyph}
\newcommand{\rcd}{$\RCstar D$\Hyph}
\newcommand{\crd}{$\CRstar D$\Hyph}

\newcommand\sT{$\star T$\Hyph}%
\newcommand\Ts{$T\star$\Hyph}%

\renewcommand{\uppercasenonmath}[1]{}

\makeatletter
\newcommand\@dotsep{4.5}
\def\@tocline#1#2#3#4#5#6#7
{\relax
		\par \addpenalty\@secpenalty\addvspace{#2}%
		\begingroup \hyphenpenalty\@M
		\@ifempty{#4}{%
			\@tempdima\csname r@tocindent\number#1\endcsname\relax
		}{%
			\@tempdima#4\relax
		}%
		\parindent\z@ \leftskip#3\relax \advance\leftskip\@tempdima\relax
		\rightskip\@pnumwidth plus1em \parfillskip-\@pnumwidth
		#5\leavevmode\hskip-\@tempdima #6\relax
		\leaders\hbox{$\m@th
		\mkern \@dotsep mu\hbox{.}\mkern \@dotsep mu$}\hfill
		\hbox to\@pnumwidth{\@tocpagenum{#7}}\par
		\nobreak
		\endgroup
}
\makeatother 

\ifx\PrintBook\undefined
	\usepackage{fancyhdr}
	\makeatletter
	\renewcommand{\@indextitlestyle}{%
		\twocolumn[\section{\indexname}]%
		\def\IndexSpace{off}%
	}
	\makeatother 
	\thanks{\href{mailto:Aleks\_Kleyn@MailAPS.org}{Aleks\_Kleyn@MailAPS.org}}
\else
	\makeatletter
	\renewcommand{\@indextitlestyle}{%
		\twocolumn[\chapter{\indexname}]%
		\def\IndexSpace{off}%
		\let\@secnumber\@empty
		\chaptermark{\indexname}%
	}
	\makeatother 
	\email{\href{mailto:Aleks\_Kleyn@MailAPS.org}{Aleks\_Kleyn@MailAPS.org}}
\fi

\ifx\SelectlEnglish\undefined
	\ifx\UseRussian\undefined
		\def\SelectlEnglish{}
	\fi
\fi

\ifx\SelectlEnglish\undefined
	\newcommand\CurrentLanguage{Russian.}%
	\author{Александр Клейн}
	\newtheorem{theorem}{Теорема}[section]
	\newtheorem{corollary}[theorem]{Следствие}
	\theoremstyle{definition}
	\newtheorem{definition}[theorem]{Определение}
	\newtheorem{example}[theorem]{Пример}
	\newtheorem{xca}[theorem]{Exercise}
	\theoremstyle{remark}
	\newtheorem{remark}[theorem]{Замечание}
	
	\newcommand\Gbasis{$G$\Hyph базис}
	\newcommand\Gcoords{$G$\Hyph координат}
	\newcommand\Gspace{$G$\Hyph пространств}
	\newcommand\xRefDef[2]
		{
		\externaldocument[#1-Russian-]{#1.Russian}[http://arxiv.org/PS_cache/#2.pdf]
		\NameDef{xRefDef#1}{}%
		}
	\makeatletter
	\newcommand\xRef[2]%
	{%
		\@ifundefined{xRefDef#1}{%
		\ref{#2}%
		}{%
		\citeBib{#1}-\ref{#1-Russian-#2}%
		}%
	}%
	\newcommand\xEqRef[2]%
	{%
		\@ifundefined{xRefDef#1}{%
		\eqref{#2}%
		}{%
		\citeBib{#1}-\eqref{#1-Russian-#2}%
		}%
	}%
	\makeatother
	\ifx\PrintBook\undefined
		\newcommand{\BibTitle}{%
			\section{Список литературы}%
		}
	\else
		\newcommand{\BibTitle}{%
			\chapter{Список литературы}%
		}
	\fi
\else
	\newcommand\CurrentLanguage{English.}%
	\author{Aleks Kleyn}
	\newtheorem{theorem}{Theorem}[section]
	
	\theoremstyle{definition}
	\newtheorem{definition}[theorem]{Definition}
	\newtheorem{example}[theorem]{Example}
	
	\theoremstyle{remark}
	\newtheorem{remark}[theorem]{Remark}
	\newcommand\Gbasis{$G$\Hyph basis}
	\newcommand\Gcoords{$G$\Hyph coordinates}
	\newcommand\Gspace{$G$\Hyph space}
	\newcommand\xRefDef[2]
		{
		\externaldocument[#1-English-]{#1.English}[http://arxiv.org/PS_cache/#2.pdf]
		\NameDef{xRefDef#1}{}%
		}
	\makeatletter
	\newcommand\xRef[2]%
	{%
		\@ifundefined{xRefDef#1}{%
		\ref{#2}%
		}{%
		\citeBib{#1}-\ref{#1-English-#2}%
		}%
	}%
	\newcommand\xEqRef[2]%
	{%
		\@ifundefined{xRefDef#1}{%
		\eqref{#2}%
		}{%
		\citeBib{#1}-\eqref{#1-English-#2}%
		}%
	}%
	\makeatother
	\ifx\PrintBook\undefined
		\newcommand{\BibTitle}{%
			\section{References}%
		}
	\else
		\newcommand{\BibTitle}{%
			\chapter{References}%
		}
	\fi
\fi

\ifx\PrintBook\undefined
\else
	\numberwithin{section}{chapter}
\fi

\numberwithin{equation}{section}
\numberwithin{figure}{section}
\numberwithin{table}{section}
\numberwithin{Item}{section}
\numberwithin{Hfootnote}{section}

\makeatletter
\newcommand\org@maketitle{}
\let\org@maketitle\maketitle
\def\maketitle{%
	\hypersetup{pdftitle={\@title}}%
	\hypersetup{pdfauthor={\authors}}%
	\hypersetup{pdfsubject=\@keywords}%
	\org@maketitle
}
\def\make@stripped@name#1{%
	\begingroup
		\escapechar\m@ne
		\global\let\newname\@empty
		\protected@edef\Hy@tempa{\CurrentLanguage #1}%
		\edef\@tempb{%
			\noexpand\@tfor\noexpand\Hy@tempa:=%
			\expandafter\strip@prefix\meaning\Hy@tempa
		}%
		\@tempb\do{%
			\if\Hy@tempa\else
				\if\Hy@tempa\else
					\xdef\newname{\newname\Hy@tempa}%
				\fi
			\fi
		}%
	\endgroup
}%
\newenvironment{enumBib}{%
	\BibTitle
	\advance\@enumdepth \@ne
	\edef\@enumctr{enum\romannumeral\the\@enumdepth}\list
	{\csname biblabel\@enumctr\endcsname}{\usecounter
	{\@enumctr}\def\makelabel##1{\hss\llap{\upshape##1}}}
}{%
	\endlist
}

\def\Chapters#1{\ChapterList#1,LastChapter,}%
\def\LastChapter{LastChapter}%
\def\ChapterList#1,{\def\temp{#1}%
	\ifx\temp\LastChapter
	\else
		\@ifundefined{#1}{%
		}{%
			\def\Semafor{on}
		}
		\expandafter\ChapterList
	\fi
}%
\newcommand{\BiblioItem}[3]
{
	\def\Semafor{off}
	\Chapters{#1}
	\ifx\Semafor\ValueOn
		\ifx\IndexState\ValueOff
			\begin{enumBib}
			\def\IndexState{on}
		\fi
		\item \label{bibitem: #2}#3%
	\fi
}
\newcommand{\OpenBiblio}
{
	\def\IndexState{off}
}
\newcommand{\CloseBiblio}
{
	\ifx\IndexState\ValueOn
		\end{enumBib}
		\def\IndexState{off}
	\fi
}

\def\StartCite{[}%
\def\citeBib#1{\showCiteBib#1,endCite,}%
\def\endCite{endCite}%
\def\showCiteBib#1,{\def\temp{#1}%
\ifx\temp\endCite
]%
\def\StartCite{[}%
\else
	\StartCite\ref{bibitem: #1}%
	\def\StartCite{, }%
\expandafter\showCiteBib%
\fi}%

\makeatother 
\newcommand{\arp}{\ar @{-->}}

\newcommand{\bundle}[4]%
{%
	\def\tempa{}%
	\def\tempb{#3}%
	\def\tempc{#1}%
	\ifx\tempa\tempb%
		\ifx\tempa\tempc%
			#2%
		\else%
			\xymatrix{#2:#1\arp[r]&#4}%
		\fi%
	\else%
		\ifx\tempa\tempc%
			#2[#3]%
		\else%
			\xymatrix{#2[#3]:#1\arp[r]&#4}%
		\fi%
	\fi%
}%
\newcommand{\AddIndex}[2]%
{%
	{\bf #1}%
	\label{index: #2}%
}%
\newcommand{\Index}[3]%
{%
	\def\Semafor{off}%
	\Chapters{#1}%
	\ifx\Semafor\ValueOn%
		\def\tempa{}%
		\def\tempb{#3}%
		\ifx\IndexState\ValueOff%
			\begin{theindex}%
			\def\IndexState{on}%
		\fi%
		\ifx\IndexSpace\ValueOn%
			\indexspace%
			\def\IndexSpace{off}%
		\fi%
		\item #2%
		\ifx\tempa\tempb%
		\else%
			\ \pageref{index: #3}%
		\fi%
	\fi%
}%
\newcommand{\SubIndex}[3]
{
	\def\Semafor{off}
	\Chapters{#1}
	\ifx\Semafor\ValueOn
		\subitem #2 \pageref{index: #3}
	\fi
}%

\makeatletter
\newcommand{\Symb}[3]
{
	\def\Semafor{off}
	\Chapters{#1}
	\ifx\Semafor\ValueOn
		\ifx\IndexState\ValueOff
			\begin{theindex}
			\def\IndexState{on}
		\fi
		\ifx\IndexSpace\ValueOn
			\indexspace
			\def\IndexSpace{off}
		\fi
		\item $\@nameuse{ViewSymbol#3}$\ \ #2
		\@nameuse{RefSymbol#3}%
	\fi
}
\def\CiteBibNote{\footnote[0]}
\makeatother

\newcommand{\SetIndexSpace}%
{%
	\def\IndexSpace{on}%
}%
\def\ValueOff{off}
\def\ValueOn{on}

\newcommand{\OpenIndex}
{
	\def\IndexState{off}
}
\newcommand{\CloseIndex}
{
	\ifx\IndexState\ValueOn
		\end{theindex}
		\def\IndexState{off}
	\fi
}
\newcommand\epigraph[2]
{
	\begin{flushright}
	\begin{minipage}{200pt}
		\Memos{#1}\CiteBibNote{#2}
	\end{minipage}
	\end{flushright}
}
\def\Memos#1{\MemoList#1//LastMemo//}%
\def\LastMemo{LastMemo}%
\def\MemoList#1//{\def\temp{#1}%
	\ifx\temp\LastMemo
	\else%
		\par\setlength{\parindent}{12pt}\textcolor{blue}{#1}%
		\expandafter\MemoList%
	\fi%
}%

\def\texPrefaceRelation{}

\xRefDef{0701.238}{math/pdf/0701/0701238}
\xRefDef{0702.561}{math/pdf/0702/0702561}
\begin{document}
\title{Fibered Correspondence}

\pdfbookmark[1]{Fibered Correspondence}{TitleEnglish}
\begin{abstract}
Base of fibered correspondence is arbitrary
correspondence. Fibered correspondence is interesting when we consider
relationship between different bundles. However composition
of fibered correspondences may not always be defined. Reduced
fibered correspondence is defined only between
fibers over the same point of base.
Reduced
fibered correspondence in bundle is called $2$\Hyph ary fibered relation.
We considered fibered equivalence and isomorphism theorem in
case of fibered morphisms.
\end{abstract}

\maketitle
\def\texPreface_Algebra{}
\ifx\texPrefaceRelation\Defined
\ifx\PrintBook\Defined
				\chapter{Preface}
\fi
\epigraph{24 - 25 August 1883
//
Pushkin once said in the circle of his friends: "Imagine what my Tat'iana has done
- she's got married. I should never have expected that of her." I could say just
the same about Anna Karenina. My characters sometimes do things that I would not wish.
In general they do what is ordinarily done in actual life and not what I want.
//
- Reminiscences of V. I. Alekseev. The date is that given by G. A. Rusanov, who
records the same saying. According to Alekseev it refers to Anna's suicide.
Tat'iana is the heroine of Pushkin's "novel in verse," Eugenii Onegin.}
{\citeBib{Tolstoi about Anna Karenina}, p 51}
\fi

\ifx\texPrefaceAlgebra\Defined
Theory of representation of algebra has long and extensive history.
During XX century representation theory became an integral part of different applications.
Transition from algebra to algebra bundle opens new opportunities.
I have ventured to write this
\ifx\PrintBook\Defined
book
\else
paper
\fi
where I want to discover new properties of algebra bundle.
\fi

\ifx\texPrefaceRelation\Defined
This story began many years ago.
When I was 15 years old I got some money for minor expenses.
I started buying math and physics books. When
I entered university I collected library of favorite books.
Among these books there was book "Universal Algebra" by Cohn
(\citeBib{Cohn: Universal Algebra})
that miraculously remained unsold. I was not thinking too hard that
it might be too early for me to read this. Reading this book
I had fallen in love with algebra.
However I dedicated my life and research to geometry 
on the edge of geometry and physics.

After more then 30 years I suddenly returned to this book.
I became curious about
definition of universal algebra in fiber of bundle.
When I started to write paper \citeBib{0702.561} I could swear that
something similar I had read when I was young. I found the book
that, as I supposed, was origin. However this book was dedicated only
to vector bundles. When I started to write this paper I finally
realized that I did not read something similar. I could not leave such facts unnoticed.

It turns out that I thought about this when I was young. However why did I wait so long time?
Why I initiated this research now? I will never get an answer to the
first question. However the answer to the second question is very simple.
For a long time I studied reference frame in general relativity.
My interest was not only classical case, but possible deviations of geometry as well.
Statements obtained in \citeBib{0702.561} show  that
in the field of my research there is a lot of unanswered questions.

I supposed to dedicate this paper to
fibered equivalence relation, because I have certain interest to it
in the future.
However, the need for clear statements involved new definitions.
Then events became unpredictable.
I expected the paper to be extremely concise and
written for a month. However this paper severely takes all my time,
changes the title and direction of the research.\footnote{Thrill
of hunt is one of the strongest passions of mankind.
I catch myself that I keep solving problems which are new for me.}
The paper dedicated to fibered equivalence relation turns into paper
dedicated to fibered binary relations.

Fibered relation is one of the most complicated subjects in the theory of fibered
algebra. Since an operation is a map,
we extend unambiguously the definition of the operation to bundle
and its sections and demand that operation is smooth.
Relation is the subset of Cartesian product.
Assuming that we defined relation only in the fiber, we are losing relationship between fibers.

I decided to repeat the procedure to determine a relation in
universal algebra.
When I started to study fibered correspondence, I realized, that I need
to change definitions in \citeBib{0702.561}.
The definition
\xRef{0702.561}{definition: reduced Cartesian product of bundle}
generates too narrow framework
to define a fibered correspondence.
May be sometimes it is enough to define correspondence only in
fiber\footnote{In particular, we define a fibered relation
introducing relation in a fiber}, however we losing
fibered morphisms.
The analysis of this situation exposes the myth
with which I have comfortably lived for all those years.
More exactly, all this time I have been trying to work out what
the base of fibered map is like. Is it an injection or an arbitrary map?
No one definition
gives clear answer on this question.
For the simplicity of perception I supposed that
base of fibered map is
injection. Actually, since we do not determine type of map of a base,
this map may be arbitrary. This leads to more wide
definition (definition
\xRef{0702.561}{definition: Cartesian product of bundle})
of Cartesian product of bundles.
On the other hand, the definition
\xRef{0702.561}{definition: Cartesian product of bundle}
presents problems to determine fibered algebra. This brings to necessity
to use two definitions of Cartesian product of bundle.
Similar considerations bring to two definitions of fibered correspondence.

Finally I return to equivalence. However,
my efforts were not wasted. My view on problem changed.

From definition \ref{definition: Fibered Correspondence}
it follows that notion of continuity is important in
definition of fibered correspondence. After this definition
I explain what does mean continuity of correspondence.
\fi


\ifx\texPrefaceAlgebra\Defined
Since a cut of the bundle
may be not defined on the whole bundle, all statements assume
a specific domain.
Such statements are modelled on statement
about existance of trivial tangent bundle on manifold.

However there exists other group of statements restricting the
domain of fibered
$\mathcal{F}$\Hyph algebra. I give appropriate examples in the text.
The explanation follows.

We suppose that transormation of fiber caused by parallel transfer
is one-to-one map.
Continuously moving along base,
we continuously move from one fiber to other.
Assumption that map between fibers is homeomorphism
waranties continuous deformation of fiber.

If we choose the point in fiber, then the trajectory
of its movement when its projection moves along
base is parallel to base.
In differential geometry such lines are called horisontal.
We also will use this definition.

Since we assume that there is structure of
$\mathcal{F}$\Hyph algebra on fiber, then we expect that
corresponding map is isomorphism
of $\mathcal{F}$\Hyph algebra.
Continuity allows to save considered structures when we use
parallel transfer, it allows to make smooth transfer from fiber to fiber.

This picture works well in the small. When we consider
finite intervals on base,
continuity becomes responsible for
impossibility to extend the structure of $\mathcal{F}$\Hyph algebra
as far as we please.
For instance, there may apear points where homeomorphism is
broken. It happens when horizontal lines have intersection or
topological propertyes of fiber change. Corresponding fiber is
called degenerate, and its projection is called point of degeneracy.

It is not easy to say how many points of degeneracy there are.
It is clear by intuition that this set is small in comparison with the base.
However, this set may prove to be essential
for the study of geometry of bundles or physical processes
associated with this bundle.

The problem to extend the solution of differential equation
is one of such events in the theory of differential equations.
At the same time, there exist two types of solution of
differential equation.
Regular solution belongs to family of functions
dependent on arbitrary constants.
Singular solution is envelope of family of regular solutions.

The problem of describing fibers of bundle, regardless whether
they are degenerate or not, has an interesting solution. Any path on the base of bundle
is map of interval $I=[0,1]$ into base of bundle.
Let us assume that fibers of bundle are not homeomorphic, but homotopic.

The holonomy group of bundle also has constrains
for structure of fibered $\mathcal F$\Hyph algebra
$\mathcal A=\bundle{}pA{}$. It is natural to assume that using
parallel transfer we have homomorphism of
$\mathcal F$\Hyph algebra from one fiber into another.
Therefore, we assume that transformation caused by
parallel transfer along loop
is homomorfism of $\mathcal F$\Hyph algebra. Thus,
everything is fine
when the holonomy group of bundle $\mathcal A$ is
subgroup of group of homomorphisms of $\mathcal F$\Hyph algebra
$A$. In this case
fibered $\mathcal F$\Hyph algebra $\mathcal A$ is called
holonomic.
Otherwise
fibered $\mathcal F$\Hyph algebra $\mathcal A$ is called
anholonomic.

From theory of vector bundles we know that there
exist fibered $\mathcal F$\Hyph algebra which is
not holonomic. At the same time, the theory of vector bundles
has answer how we can work with
anholonomic fibered $\mathcal F$\Hyph algebra.

We apply this remark also to the theory of representations of
fibered $\mathcal F$\Hyph algebra.

A new design is illustrated through the use of corresponding diagrams.
\fi

Where it is possible
I use the same notation for operations and relations
as we use them in the set theory. It does not bring to
ambiguity because we use different notation for set and
bundle. I use the same letter in different alphabets to denote
bundle and fiber.

We assume that projection of bundle, section and fibered map are smooth maps.

In mathematical literature there are two customs to write product of
mappings and correspondences. Some authors write product of
mappings in the same order as arrows follow on diagram.
While others prefer to write mappings in opposite
order. When we read papers and books the first thing what we
need is to put attention what order of factors is used by author.

The case is clearer when author writes an action of mapping
over set. In this case the author writes set and mapping in such
order that, when we write brackets, we gets right order.
For instance, let us consider diagram
\[
\xymatrix{
A\ar[r]^f&B\ar[r]^g&C
}
\]
and let $D\subset A$. Since we write product of mappings
as $fg$, then image of set $D$ has form $Dfg=(Df)g$.
Since we write product of mappings
as $gf$, then image of set $D$ has form $gfD=g(fD)$.

This is my start point. Based on convention
from remark \xRef{0701.238}{remark: left and right matrix notation},
I will assume oportunity to read expression from right to left
and from left to right.
\def\texBundleRelation{}
\ifx\PrintBook\Defined
				\chapter{Fibered Relation}
\fi

			\section{Correspondence}

		\begin{definition}
Let
\[
\xymatrix{
A\ar[d]^\Phi\ar[rr]^\Psi&&B\ar[d]^\Sigma\\
C\ar[rr]^\Theta&&D
}
\]
be diagram arrows of which represent correspondences.
\AddIndex{Diagram of correspondences}{diagram of correspondences} is called
\AddIndex{commutative}{commutative diagram of correspondences} when image of any subset
of set $A$ in set $D$ does not depend on way in diagram.
		\qed
		\end{definition}

The definition of category is not specific on the question whether morphism is a map.
We can study category, sets of which are objects and
correspondences from one set into another are morphisms.

		\begin{definition}
Let $\Phi$ be a correspondence from a set $A$ to a set $B$.
Let $C\subseteq A$.
The correspondence
\symb{\Phi/C}0{restriction of correspondence}
\[
\ShowSymbol{restriction of correspondence}=\{(c,b)\in\Phi:c\in C\}
\]
is called a \AddIndex{restriction of the correspondence $\Phi$ to the set $C$}
{restriction of correspondence}.\footnote{We
make this definition similar to the definition from \citeBib{Bourbaki: Set Theory}, p. 82}
correspondence $\Phi$ is called an \AddIndex{extension of correspondence}
{extension of correspondence} $\Phi/C$.
		\qed
		\end{definition}

		\begin{definition}
		\label{definition: correspondence continuous on the set}
Let sets $A$ and $B$ be topological spaces.

Correspondence $\Phi$ from set $A$ to a set $B$
is said to be
\AddIndex{continuous on the set}{correspondence continuous on the set}
$C\subset A$,
if, given open set $V$,
$\Phi C\subset V\subset B$, there is an open set
$U$, $C\subset U\subset A$,
such that $\Phi U\subseteq V$.

Correspondence $\Phi$ from set $A$ to a set $B$
is said to be \AddIndex{continuous}{continuous correspondence},
if, given open set
$V\subset B$, there is an open set $U\subset A$
such that $\Phi U\subseteq V$.
		\qed
		\end{definition}

Following after \citeBib{Bourbaki: General Topology 1},
I define continuity based on definition of
limit of filter. Since an image under correspondence is not
a point, but a set, I little changed definitions and theorems.

		\begin{definition}
		\label{definition: filter converges}
Let $X$ be a topological space and $\mathfrak{F}$ a
filter on $X$. A set $A\subset X$ is said to be a
\AddIndex{limit set}{limit set of filter} or
\AddIndex{limit of filter}{limit of filter}
$\mathfrak{F}$, if $\mathfrak{F}$ is finer then the
neighborhood filter $\mathfrak{B}(A)$ of set $A$.
Filter $\mathfrak{F}$ is also said to
\AddIndex{converge}{filter converges} to $A$,
\symb{\mathfrak{F}\rightarrow A}1{filter converges}.

The set $A$ is said to be a limit of filter base $\mathfrak{B}$
on $X$ and $\mathfrak{B}$ is said to converge to $A$,
if the filter whose base is $\mathfrak{B}$ converges to $A$.
		\qed
		\end{definition}

		\begin{theorem}
		\label{theorem: filter base converges to set}
A filter base $\mathfrak{B}$ on topological space
$X$ converges to set $A\subset X$ iff
every set of fundamental system of
neighborhoods of set $X$ contains a set of $\mathfrak{B}$.
		\end{theorem}
		\begin{proof} 
If a filter $\mathfrak{F}$ converges to set $A$, then every filter
finer than $\mathfrak{F}$ also converges to $A$,
by definition \ref{definition: filter converges}.
Let $\Phi$ be a set of filters on $X$, all of which converges to
set $A$. The neighborhood filter $\mathfrak{B}(A)$
is coarser than all filters of $\Phi$, hence $\mathfrak{B}(A)$
is coarser than their intersection $\mathfrak{G}$. Therefore,
$\mathfrak{G}$ converges to set $A$.
		\end{proof}

		\begin{definition}
		\label{definition: limit of correspondence with respect to the filter}
Let $\Phi$ be a correspondence from set $X$ to a
topological space $Y$. Let $\mathfrak{F}\rightarrow A$ be
a filter on $X$. The set $B\subset Y$ is said to be
\AddIndex{limit of correspondence with respect to the filter}
{limit of correspondence with respect to the filter}
$\mathfrak{F}$
\symb{\lim_{\mathfrak{F}\rightarrow A}\Phi(\mathfrak{F})}
0{limit of correspondence with respect to the filter}
\[
\ShowSymbol{limit of correspondence with respect to the filter}=B
\]
if $B$ is limit of the filter base
$\Phi(\mathfrak{F})$.
		\qed
		\end{definition}

		\begin{theorem}
		\label{theorem: limit of correspondence}
A set $B\subset Y$
is a limit of correspondence $\Phi$ with respect to the filter
$\mathfrak{F}$ iff for each neighborhood $V$ of set $B$
in $Y$, there is a set $M\in\mathfrak{F}$
such that $\Phi M\subset V$.
		\end{theorem}
		\begin{proof} 
The statement is corollary of definition
\ref{definition: limit of correspondence with respect to the filter}
and theorem \ref{theorem: filter base converges to set}.
		\end{proof}

		\begin{theorem}
A correspondence $\Phi$ from
a topological space $X$ to
a topological space $Y$ is continuous
on the set $A\subset X$ iff
\[
\lim_{\mathfrak{F}\rightarrow A}\Phi(\mathfrak{F})=\Phi A
\]
		\end{theorem}
		\begin{proof} 
The statement is corollary of definition
\ref{definition: correspondence continuous on the set}
and theorem \ref{theorem: limit of correspondence}.
		\end{proof}

Let $\Phi$ be continuous correspondence from
topological space $X$ to
topological space $Y$. Let
$(a,b)\in\Phi$. Suppose $V\subset Y$ is an open set,
$\Phi\{a\}\subset V$. In particular, $b\in V$.
According to definition
\ref{definition: correspondence continuous on the set}
there exists an open set $U\subset X$, $a\in U$,
$\Phi U\subset V$. Therefore, there exist $a'\in U$,
$b'\in V$, $(a',b')\in\Phi$.

			\section{Fibered Correspondence}

		\begin{definition}
		\label{definition: fibered subset}
Let $\bundle{\mathcal{A}}aAN$ and $\bundle{\mathcal{B}}bBM$ be bundles.
Suppose bundle map is defined by diagram
	\[
\xymatrix{
\mathcal{A}\arp[d]^{\bundle{}aA{}}\ar[rr]^f
&&\mathcal{B}\arp[d]^{\bundle{}{b}{B}{}}\\
N\ar[rr]^F&&M
}
	\]
where maps $f$ and $F$ are injections. Then bundle $\bundle{}aA{}$ is called
\AddIndex{fibered subset}{fibered subset} or \AddIndex{subbundle}{subbundle} of $\bundle{}bB{}$.
We also use notation
\symb{\bundle{}{a}{A}{}\subseteq\bundle{}{b}{B}{}}1{fibered subset} or
\symb{\mathcal{A}\subseteq\mathcal{B}}1{subbundle}.

Suppose we defined the $\mathcal{F}$\Hyph algebra
on sets $A$ and $B$ and map $f$
is homomorphism of fibered algebra. Then the bundle $\bundle{}aA{}$
is called fibered subalgebra of fibered algebra $\bundle{}bB{}$.
		\qed
		\end{definition}

Without loss of generality we assume that $A\subseteq B$, $N\subseteq M$.

		\begin{definition}
		\label{definition: Fibered Correspondence}
Let $\bundle{\mathcal{A}}aAM$ and $\bundle{\mathcal{B}}bBN$ be bundles.
Fibered subset $\bundle{\mathcal{F}}fF\Phi$ of bundle $\mathcal{A}\times\mathcal{B}$
is called \AddIndex{fibered correspondence from $\mathcal{A}$ to $\mathcal{B}$}{fibered correspondence from A to B}.

If $\mathcal{A}=\mathcal{B}$, then fibered correspondence $\mathcal{F}$ is called
\AddIndex{fibered correspondence in $\mathcal{A}$}{fibered correspondence in A}.
		\qed
		\end{definition}

According to definitions \ref{definition: fibered subset} and
\ref{definition: Fibered Correspondence}
we can represent fibered correspondence using the diagram
\begin{equation}
\label{eq: Fibered Correspondence, definition}
\xymatrix{
\mathcal{F}\ar[rr]^j\arp[dd]^{\bundle{}fF{}}&&
\mathcal{A}\times\mathcal{B}
\arp@/^1pc/[dd]^{\bundle{}bB{}}
\arp@/_1pc/[dd]_{\bundle{}aA{}}\\
&&\times\\
\Phi\ar[rr]^i&&M\times N
}
\end{equation}
where $i$ and $j$ are continuous injections.
We assume that a set $U\subset\Phi$ is open
iff there exists an open set
$V\subset M\times N$ such that $U=V\cap\Phi$.
We assume that a set $\mathcal{U}\subset\mathcal{F}$ is open
iff there exists an open set
$\mathcal{V}\subset\mathcal{A}\times\mathcal{B}$
such that $\mathcal{U}=\mathcal{V}\cap\mathcal{F}$.

Correspondence $\Phi$ is called
\AddIndex{base of fibered correspondence}{base of fibered correspondence}
$\mathcal{F}$, and fibered correspondence $\mathcal{F}$
is called \AddIndex{lift of correspondence}{lift of correspondence} $\Phi$.

We use diagrams of fibered correspondences as
we use diagrams of maps and correspondences.
In this case, in the diagram
we can additionally show a projection on the base.
Thus we can represent fibered correspondence using diagram
\begin{equation}
\label{eq: Fibered Correspondence, representation}
\xymatrix{
\mathcal{A}\ar[rr]^{\bundle{}pF{}}="AR1"\arp[d]^a&&\mathcal{B}\arp[d]^b\\
M\ar[rr]_f="AR2"&&N\\
\arp^{\bundle{}pF{}} "AR1";"AR2"
}
\end{equation}%
Choice of diagram \eqref{eq: Fibered Correspondence, definition} or
\eqref{eq: Fibered Correspondence, representation}
depends on problem which we consider.

To study fibered correspondence from $\mathcal{A}$ to $\mathcal{B}$
I use maps of the manifold $M$ and $N$ in which both bundles are
trivial. This will make it possible to review fibered correspondence in detail 
without loss of generality.


According to definition \xRef{0702.561}{definition: Cartesian product of bundle},
given $x\in M$, $y\in N$,
I represent points of a fiber $(A\times B)_{(x,y)}$
of the bundle $\mathcal{A}\times\mathcal{B}$ as
tuple $(x,y,p,q)$, $p\in A_x$, $q\in B_y$. According to our assumption,
$\Phi\subseteq M\times N$, $F\subseteq A\times B$.
Therefore, we can consider $\Phi$ as correspondence from $M$ to $N$
and $F$ as correspondence from $A$ to $B$. In particular, $F_{(x,y)}\subseteq A_x\times B_y$.
The point $(x,p)\in\mathcal{A}$ is in correspondence $\mathcal{F}$ with point
$(y,q)\in\mathcal{B}$, if point $x\in M$ is in correspondence $\Phi$ with
point $y\in N$ and point $p\in A_x$ is in correspondence $F_{(x,y)}$ with
point $q\in B_y$.

Correspondence in fiber depends on selected fiber. For instance, let
$\bundle{\mathcal{A}}aR{R^2}$ and $\bundle{\mathcal{B}}bRR$ be bundles.
Let $(x,y)\in M=R^2$ and $z\in N=R$. Assume that point
$(x,y,p)\in\mathcal{A}$ is in correspondence $\mathcal{F}$
with point $(z,q)\in\mathcal{B}$, if we can represent $q$ as
\[
q=(z^2+x^2+y^2+p^2+1)n
\]
where $n$ is an arbitrary integer. The relation is different
in different fibers, however we can define bijection
between $F$ and $F_{(x,y,z)}$
for arbitrary fiber.\footnote{One can easily observe that small change in
correspondence may bring to fact that in some fibers
correspondence will be singular.}

Let us consider set $\Gamma(\mathcal{F})$. Since we choose
clutching functions of bundles $\mathcal{A}$ and $\mathcal{B}$,
we may represent an element $\Gamma(\mathcal{F})$ as
$(x,y,p(x),q(y))$. Little permutation
leads to the record $((x,p(x)),(y,q(y)))$. This leads one to assume that
we wrote correspondence from set $\Gamma(\mathcal{A})$
to set $\Gamma(\mathcal{B})$. However this is not so.
Let us choose value of $x$. Then when we change $y$, tuple $(y,q(y))$
describes section $q\in\Gamma(\mathcal{B})$. Since this section
depends on choice of tuple $(x,p(x))$, we cannot establish
correspondence from set $\Gamma(\mathcal{A})$
to set $\Gamma(\mathcal{B})$.

Fibered correspondence
loses some important properties of correspondence. For instance,
if $C\subseteq A$ and $\Phi$ is correspondence from $A$ to $B$,
then we can define the image of the set $C$ under correspondence $\Phi$
using law
\[
\Phi C=\{b\in B:(a,b)\in\Phi,a\in C\}
\]
However in the case of a fibered correspondence the image of a bundle
is not necessarily a bundle, which is due to the fact
that there exists a possibility that $(x_1,y)\in\Phi$, $(x_2,y)\in\Phi$.
At the same time, generally speaking,
$$\Phi_{(x_1,y)}A_{x_1}\ne \Phi_{(x_2,y)}A_{x_2}$$ even
$$\Phi_{(x_1,y)}A_{x_1}\subseteq B_y$$ $$\Phi_{(x_2,y)}A_{x_2}\subseteq B_y$$
Corollary of this is impossibility to define
a composition of fibered correspondences in general case.
Similar statement holds for fibered morphisms when base
of morphism is not injection.

		\begin{theorem}
		\label{theorem: image of fibered correspondence}
Let us consider fibered correspondence $\mathcal{F}$
\[
\xymatrix{
\mathcal{A}\arp[d]^{\bundle{}pA{}}\ar[rr]^{\mathcal{F}}="AR1"
&&\mathcal{B}\arp[d]^{\bundle{}qB{}}\\
M\ar[rr]_f="AR2"&&N
\arp^{\bundle{}sF{}} "AR1";"AR2"
}
\]
from bundle $\mathcal{A}$ to bundle $\mathcal{B}$,
base $f$ of which is injection.
Let the bundle
\[\bundle{\mathcal{C}}aCL\]
is subbundle of bundle
$\mathcal{A}$. We define the image
of the bundle $\mathcal{C}$ under fibered correspondence $\mathcal{F}$
according to law
\begin{align*}
\mathcal{F}\mathcal{C}=\{(y,b):&y\in N,\exists x\in M,y=f(x),\\
&b\in B_y,\exists a\in A_x,(a,b)\in F_{(x,y)}\}
\end{align*}
Image of bundle $\mathcal{C}$ under fibered
correspondence $\mathcal{F}$
is subbundle of bundle $\mathcal{B}$.
		\end{theorem}
		\begin{proof}
To prove the statement, we need to show that all
$D_y=F_{(x,y)}C_x$ are homeomorphic.

Let us consider the following diagram
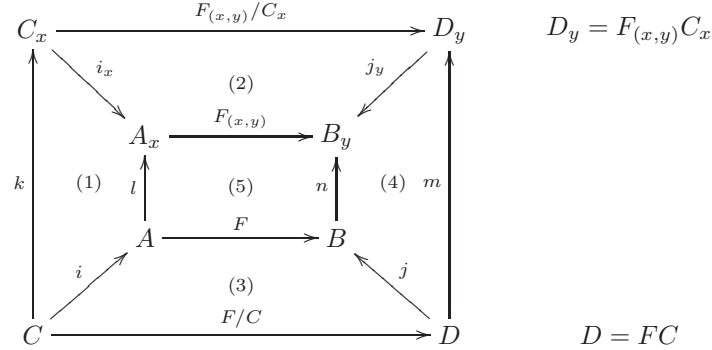
\[
\xymatrix{
C_x\ar[dr]^{i_x}\ar[rrrr]^{F_{(x,y)}/C_x}\ar@{}[drrrr]|{(2)}&&&\ar@{}[dddr]|{(4)}&
D_y\ar[dl]_{j_y}&D_y=F_{(x,y)}C_x\\
&A_x\ar[rr]^{F_{(x,y)}}&&B_y&\\
&A\ar[rr]^F\ar[u]^l\ar@{}[urr]|{(5)}&&B\ar[u]^n&\\
C\ar[ur]^i\ar[uuu]^k\ar[rrrr]^{F/C}\ar@{}[uuur]|{(1)}\ar@{}[urrrr]|{(3)}&&&&
D\ar[ul]_j\ar[uuu]^m&D=FC
}
\]
$i$, $i_x$, $j$, and $j_y$ are injections, $k$, $l$, and
$n$ are bijections.
We need to prove, that $m$ is bijection. We assume that correspondence
$F$ is nonsingular in selected fiber.

Bijection $l$ means that we can
enumerate points of set $A_x$ using points of set $A$.
Bijection $k$ means that we can
enumerate points of set $C_x$ using points of set $C$.
Injection $i$ means that $C\subseteq A$.
Injection $i_x$ means that $C_x\subseteq A_x$.
Therefore, for each point $p\in A$, point $p_x\in A_x$ is defined uniquely.
Commutativity of diagram (1) means that
	\begin{equation}
p\in C\Leftrightarrow p_x\in C_x
	\label{eq: Fibered correspondence, p_x in C_x}
	\end{equation}

Bijection $n$ means that we can
enumerate points of the set $B_x$ using points of the set $B$.
Injection $j$ means that $D\subseteq B$.
Injection $j_y$ means that $D_y\subseteq B_y$.
Therefore, for each point $q\in B$, point $q_y\in B_y$ is defined uniquely.

Bijections $l$ and $n$ and commutativity of diagram (5) means that we can
enumerate points of correspondence $F_{(x,y)}$ using points
of correspondence $F$\footnote{Requirement of nonsingularity
of correspondence in fiber is very important.
Commutativity of diagram is
broken when correspondence in a fiber is singular.}
	\begin{equation}
(p,q)\in F\Leftrightarrow (p_x,q_y)\in F_{(x,y)}
	\label{eq: Fibered correspondence, (p_x,q_y) in F_(x,y)}
	\end{equation}

According to definition, $(p,q)\in F/C$ when $p\in C$ and $(p,q)\in F$.
According to \eqref{eq: Fibered correspondence, p_x in C_x} $p_x\in C_x$.
According to \eqref{eq: Fibered correspondence, (p_x,q_y) in F_(x,y)} $(p_x,q_y)\in F_{(x,y)}$.
According to definition, $(p_x,q_y)\in F_{(x,y)}/C_x$.
Therefore, $q_y\in D_y$, and map $m$ is injection.
		\end{proof}

		\begin{theorem}
Let us consider fibered correspondence $\bundle{}sF{}$
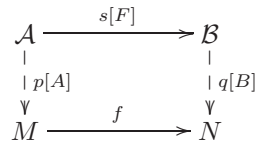
\[
\xymatrix{
\mathcal{A}\arp[d]^{\bundle{}pA{}}\ar[rr]^{\bundle{}sF{}}&&
\mathcal{B}\arp[d]^{\bundle{}qB{}}\\
M\ar[rr]^f&&N
}
\]
from bundle $\mathcal{A}$ to bundle $\mathcal{B}$
and fibered correspondence $\bundle{}tH{}$
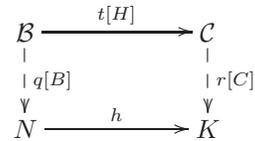
\[
\xymatrix{
\mathcal{B}\arp[d]^{\bundle{}qB{}}\ar[rr]^{\bundle{}tH{}}&&
\mathcal{C}\arp[d]^{\bundle{}rC{}}\\
N\ar[rr]^h&&K
}
\]
from bundle $\mathcal{B}$ to bundle $\mathcal{C}$.
Let bases of fibered correspondences $\bundle{}sF{}$ and $\bundle{}tH{}$
are injections. We define a
\AddIndex{composition of fibered correspondences}
{composition of fibered correspondences}\footnote{Composition of
correspondences $\Phi$ and $\Psi$
is determined even when $\Phi$ is correspondence to set $B$,
and $\Psi$ is correspondence from
set $C$. However we assume without loss of generality that
$\Phi$ is correspondence to set $B\cap C$, and $\Psi$ is correspondence from
set $B\cap C$. This allows establishing connection between
composition of correspondences and composition of fibered correspondences.}
$\mathcal{H}$ and  $\mathcal{F}$
\symb{\bundle{}tH{}\circ\bundle{}sF{}}0{composition of fibered correspondences}
\begin{align*}
\ShowSymbol{composition of fibered correspondences}
=\{(x,z,a,c):&x\in M, z\in K,
\exists y\in N,y=f(x),z=h(y),\\
&a\in A_x, c\in C_z, \exists b\in B_y, (a,b)\in F_{(x,y)},(b,c)\in H_{(y,z)}\}
\end{align*}
		\end{theorem}
		\begin{proof}
Let us consider the following diagram
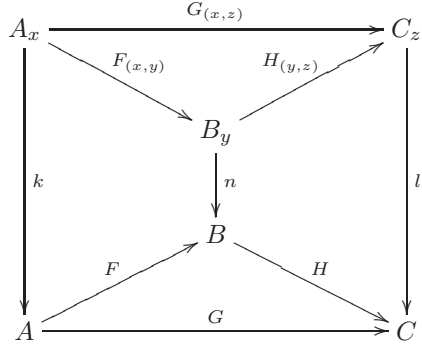
\[
\xymatrix{
A_x\ar[rrrr]^{G_{(x,z)}}\ar[drr]^{F_{(x,y)}}\ar[ddd]^k&&&&C_z\ar[ddd]^l\\
&&B_y\ar[urr]^{H_{(y,z)}}\ar[d]^n&&\\
&&B\ar[drr]^H&&\\
A\ar[rrrr]^G\ar[urr]^F&&&&C
}
\]
$k$, $l$, and $n$ are bijections.
We assume that
correspondence $F$ is nonsingular in a selected fiber.

Like in the proof of theorem
\ref{theorem: image of fibered correspondence}
commutativity of vertical diagram means that
we can
enumerate points of correspondence $F_{(x,y)}$ using points of correspondence $F$,
points of correspondence $H_{(y,z)}$ using points of correspondence $H$,
points of correspondence $G_{(x,z)}$ using points of correspondence $G$.

Commutativity of lower diagram means that
$G=H\circ F$.
Commutativity of upper diagram means that
$G_{(x,z)}=H_{(y,z)}\circ F_{(x,y)}$.
		\end{proof}

		\begin{theorem}
		\label{theorem: associative law, composition of fibered correspondences}
Let $\bundle{}sF{}$ be fibered correspondence
from bundle $\mathcal{A}$ to bundle $\mathcal{B}$,
$\bundle{}tH{}$ be fibered correspondence
from bundle $\mathcal{B}$ to bundle $\mathcal{C}$
and $\bundle{}rG{}$ be fibered correspondence
from bundle $\mathcal{C}$ to bundle $\mathcal{D}$.
If there exist compositions of fibered correspondences
	\begin{equation}
\bundle{}tH{}\circ\bundle{}sF{}
	\label{eq: composition of fibered correspondences, F H}
	\end{equation}
and
	\begin{equation}
\bundle{}rG{}\circ\bundle{}tH{}
	\label{eq: composition of fibered correspondences, H G}
	\end{equation}
then there exists compositions
$\bundle{}rG{}\circ(\bundle{}tH{}\circ\bundle{}sF{})$
and $(\bundle{}rG{}\circ\bundle{}tH{})\circ\bundle{}sF{}$.
In this case composition of fibered correspondences holds
\AddIndex{associative law}
{associative law, composition of fibered correspondences}
\[
\bundle{}rG{}\circ(\bundle{}tH{}\circ\bundle{}sF{})=
(\bundle{}rG{}\circ\bundle{}tH{})\circ\bundle{}sF{}
\]
		\end{theorem}
		\begin{proof}
Existence of composition
\eqref{eq: composition of fibered correspondences, F H} and
\eqref{eq: composition of fibered correspondences, H G} means
that base $f$ of fibered correspondence $\bundle{}sF{}$,
base $h$ of fibered correspondence $\bundle{}tH{}$ and
base $g$ of fibered correspondence $\bundle{}rG{}$
are injections. In this case there exist compositions of mappings
$h\circ f$ and $g\circ h$ which also are injections.
Therefore, there exist compositions of mappings
$g\circ(h\circ f)$ and $(g\circ h)\circ f$ which are
injections and satisfy to law
\[
g\circ(h\circ f)=(g\circ h)\circ f
\]
Therefore, there exist mapping of base of bundle $\mathcal{A}$ to
base of bundle $\mathcal{D}$, and this mapping is injection.

Existence of compositions
\eqref{eq: composition of fibered correspondences, F H} and
\eqref{eq: composition of fibered correspondences, H G} means
that there exist compositions of correspondences $H\circ F$ and
$G\circ H$. Therefore, there exist compositions of correspondences
$G\circ(H\circ F)$ and $(G\circ H)\circ F$ which hold
to the associative law. Therefore, correspondence from fiber
of the bundle $\mathcal{A}$ to fiber
of the bundle $\mathcal{D}$ is determined uniquely.
		\end{proof}

		\begin{definition}
Let $\bundle{\mathcal{F}}sFf$ be fibered correspondence
\[
\xymatrix{
\mathcal{A}\arp[d]^p\ar[rr]^{\mathcal{F}}="AR1"&&\mathcal{B}\arp[d]^q\\
M\ar[rr]_f="AR2"&&N
\arp^{\bundle{}sF{}} "AR1";"AR2"
}
\]
from bundle $\mathcal{A}$ to bundle $\mathcal{B}$
and mapping $f$ be injection.
Then there exist
\AddIndex{inverse fibered correspondence}{inverse fibered correspondence}
\symb{\mathcal{F}^{-1}}0{inverse fibered correspondence, 1}
\symb{\bundle{}s{F^{-1}}{}}0{inverse fibered correspondence, 2}
$\bundle{\ShowSymbol{inverse fibered correspondence, 1}}s{F^{-1}}{f^{-1}}$
\[
\xymatrix{
\mathcal{B}\arp[d]^p
\ar[rr]^{\ShowSymbol{inverse fibered correspondence, 1}}="AR1"
&&\mathcal{A}\arp[d]^q\\
N\ar[rr]_{f^{-1}}="AR2"&&M
\arp^{\ShowSymbol{inverse fibered correspondence, 2}} "AR1";"AR2"
}
\]
		\qed
		\end{definition}

		\begin{definition}
Let $\bundle{\mathcal{A}}pAM$ be bundle.
Fibered correspondence
\symb{\Delta_\mathcal{A}}0{diagonal in bundle, 1}
\symb{\bundle{}{r_\Delta}{\Delta_A}{}}0{diagonal in bundle, 2}
$$\bundle{\ShowSymbol{diagonal in bundle, 1}}{r_\Delta}{\Delta_A}{\Delta_M}$$
where projection $\ShowSymbol{diagonal in bundle, 2}$ is defined by law
$$r_\Delta((x,p),(x,p))=(x,x)$$ is called
\AddIndex{diagonal in bundle}{diagonal in bundle} $\mathcal{A}$.
		\qed
		\end{definition}

		\begin{theorem}
Let $\mathcal F$ be fibered correspondence
from bundle $\mathcal{A}$ to bundle $\mathcal{B}$,
$\mathcal H$ be fibered correspondence
from bundle $\mathcal{B}$ to bundle $\mathcal{C}$,
projections of fibered correspondences $\mathcal F$ and
$\mathcal H$ be injections.
The following laws are valid for
fibered correspondences $\mathcal F$ and
$\mathcal H$
\[
(\mathcal H\circ\mathcal F)^{-1}=\mathcal F^{-1}\circ\mathcal H^{-1}
\]
\[
(\mathcal F^{-1})^{-1}=\mathcal F
\]
\[
\mathcal F\circ\Delta_\mathcal{A}=\Delta_\mathcal{B\circ\mathcal F}=\mathcal F
\]
		\end{theorem}
		\begin{proof}
The proof of theorem is similar to proof of theorem
\ref{theorem: associative law, composition of fibered correspondences}.
We check each statement for base and fiber.
		\end{proof}

			\section{Fibered Correspondence of Homomorphism}

Let transition functions $f_{\alpha\beta}$ determine bundle $\mathcal{A}$
over base $M$.
Let us consider
maps $U_\alpha\in M$ and $U_\beta\in M$, $U_\alpha\cap U_\beta\ne\emptyset$.
Point $p\in\mathcal{A}$ has representation $(x,p_\alpha)$
in map $U_\alpha$ and representation
$(x,p_\beta)$ in map $U_\beta$.
Let transition functions $g_{\epsilon\delta}$ determine bundle $\mathcal{B}$
over base $N$.
Let us consider maps $V_\epsilon\in N$ and $V_\delta\in N$,
$V_\epsilon\cap V_\delta\ne\emptyset$.
Point $q\in\mathcal{B}$ has representation $(y,q_\epsilon)$
in map $V_\epsilon$ and representation
$(y,q_\delta)$ in map $V_\delta$.
Therefore,
\[
p_\alpha=f_{\alpha\beta}(p_\beta)
\]
\[
q_\epsilon=g_{\epsilon\delta}(q_\delta)
\]
When we move from map $U_\alpha$ to map $U_\beta$
and from map $V_\epsilon$ to map $V_\delta$,
representation of correspondence changes according to the law
\[
(x,y,p_\alpha,q_\epsilon)=(x,y,f_{\alpha\beta}(p_\beta),g_{\epsilon\delta}(q_\delta))
\]
This is consistent with the transformation
when we move from map $U_\alpha\times V_\epsilon$ to map $U_\beta\times V_\delta$
in the bundle $\mathcal{A}\times\mathcal{B}$.

Actually this is unusual fact.
It is enough to consider at least one operation in algebra
to demand transition functions to be homomorphisms of algebra.
However, we have here arbitrary condition, arbitrary correspondence.
And everything is OK.

As opposed to operation,
the only property which correspondence holds is belonging
of point to a certain set.
There is no reason to search for constraint for type of map
until we consider one-to-one map.
On the other hand, an operation in
the algebra can
be represented as correspondence. For instance, given
vector space,
we can consider correspondence of vector $a$ to vector $b$ such, that $a=3b$.
According to definition a linear transformation holds this correspondence.
What happens in case of
nonlinear transformation. Let us consider coordinate transformation
$a_i\rightarrow a_i^2$. Does linear relation between vectors hold? No.
Does the correspondence hold? Yes, however it is
expressed by another way. Vectors are
the same, even their coordinates change.

From this follows, that we apply constraints
to transition functions
only when we apply constraints to correspondence.

		\begin{definition}
Correspondence $\Phi$ from $\mathcal{F}$\Hyph algebra $A$ into $\mathcal{F}$\Hyph algebra $B$
is called \AddIndex{correspondence of homomorphism}
{correspondence of homomorphism}, if
\[
(\omega(a_1,...,a_n),\omega(b_1,...,b_n))\in\Phi
\]
for each n\Hyph ari operation $\omega$ and any
set of elements $a_1$, ..., $a_n\in A$, $b_1$,
..., $b_n\in B$ such that
\[
(a_1,b_1)\in\Phi, ..., (a_n,b_n)\in\Phi
\]
		\qed
		\end{definition}

		\begin{definition}
Correspondence $\Phi$ from $\mathcal{F}$\Hyph algebra $A$ into $\mathcal{F}$\Hyph algebra $B$
is called \AddIndex{fibered correspondence of homomorphism}
{fibered correspondence of homomorphism}, if
\[
(\omega(a_1,...,a_n),\omega(b_1,...,b_n))\in\Phi
\]
for each n\Hyph ari operation $\omega$ and any
set of elements $a_1$, ..., $a_n\in A$, $b_1$,
..., $b_n\in B$ such that
\[
(a_1,b_1)\in\Phi, ..., (a_n,b_n)\in\Phi
\]
		\qed
		\end{definition}

			\section{Reduced Fibered Correspondence}

		\begin{definition}
		\label{definition: reduced fibered correspondence}
Let $\bundle{\mathcal{A}}aAM$ and $\bundle{\mathcal{B}}bBM$ be bundles over $M$.
Fibered subset $\bundle{\mathcal{F}}fFN$ of bundle $\mathcal{A}\odot\mathcal{B}$
is called \AddIndex{reduced fibered correspondence from $\mathcal{A}$ to $\mathcal{B}$}
{reduced fibered correspondence from A to B}.

If $\mathcal{A}=\mathcal{B}$, then reduced fibered correspondence $\mathcal{F}$ is called
\AddIndex{reduced fibered correspondence in $\mathcal{A}$}{reduced fibered correspondence in A}.
		\qed
		\end{definition}

According to definition \ref{definition: fibered subset} and
\ref{definition: reduced fibered correspondence}
we can represent reduced fibered correspondence
using diagram
\begin{equation}
\label{eq: reduced fibered correspondence, definition}
\xymatrix{
\mathcal{F}\ar[rr]^j\arp[dd]_{\bundle{}fF{}}&&
\mathcal{A}\odot\mathcal{B}
\arp@/^1pc/[dd]^{\bundle{}bB{}}
\arp@/_1pc/[dd]_{\bundle{}aA{}}\\
&&\odot\\
M\ar[rr]^{\mathrm{id}}&&M
}
\end{equation}
where $j$ is continuous injection.
We suppose, that the set $\mathcal{U}\subset\mathcal{F}$
is open iff there exists set
$\mathcal{V}\subset\mathcal{A}\odot\mathcal{B}$
such that $\mathcal{U}=\mathcal{V}\cap\mathcal{F}$.

We define reduced fibered correspondence only for points of the same fiber.
To study reduced fibered correspondence from $\mathcal{A}$ to $\mathcal{B}$
I use maps of the manifold $M$ in which both bundles are
trivial. This will make it possible to review
reduced fibered correspondence in detail 
without loss of generality.

According to definition \xRef{0702.561}{definition: reduced Cartesian product of bundle},
I represent points of a fiber $(A\times B)_x$ of the bundle $\mathcal{A}\odot\mathcal{B}$ as
tuple $(x,p,q)$, $p\in A_x$, $q\in B_x$.
We define reduced fibered correspondence only for points of the same fiber.
The point $(x,p)\in\mathcal{A}$ is in
reduced fibered correspondence $\mathcal{F}$ with point
$(x,q)\in\mathcal{B}$, if $x\in N\subseteq M$
and point $p\in A_x$ is in correspondence $F_x$ with
point $q\in B_x$.
Therefore, we can consider $F_x$ as correspondence from $A_x$ to $B_x$.
In particular, $F_x\subseteq A_x\times B_x$.

Since reduced fibered correspondence inherits
topology of bundle $\mathcal{A}\odot\mathcal{B}$,
we can consider topology of reduced fibered correspondence.
Let point $p\in A_x$ be in correspondence $F_x$ with
point $q\in B_x$. Let
$\mathcal{V}\subset\mathcal{A}\odot\mathcal{B}$ be
an open set such that $(x,p,q)\in\mathcal{V}$.
According to \citeBib{Bourbaki: General Topology 1}, page 44,
there exist open sets $U\subset M$,
$V\subset A\times B$ such that
$x\in U$, $(p,q)\in V\cap F$.
Therefore, there exist $x'\in U$,
$(p',q')\in V\cap F$.
If $x\ne x'$, then we can express this fact
as the statement about continuous dependence of correspondence on fiber.
However
continuity of correspondence in fiber does not follow
from continuous dependence of correspondence on fiber. 

		\begin{theorem}
		\label{theorem: reduced fibered correspondence}
Given reduced fibered correspondence $\bundle{\mathcal{F}}fFN$,
we determine uniquely a fibered correspondence
$\bundle{\mathcal{F}}{f_\Delta}F{\Delta_N}$,
which is lift of diagonal $\Delta_N$,
and isomorphism of bundles
	\begin{equation}
\xymatrix{
\mathcal{F}\arp[d]^{\bundle{}{f_\Delta}F{}}\ar[rr]^{id}
&&\mathcal{F}\arp[d]^{\bundle{}fF{}}\\
\Delta_N\ar[rr]^{\pi}&&N
}
	\label{eq: Reduced Fibered correspondence}
	\end{equation}
over base $\pi:\Delta_N\rightarrow N$.
		\end{theorem}
		\begin{proof}
According to design, $F_{(x,x)}=F_x$.
		\end{proof}

We will write diagram \eqref{eq: Reduced Fibered correspondence}
in more compact form
\[
\xymatrix{
&\mathcal{F}\arp[ld]_{\bundle{}{f_\Delta}F{}}
\arp[rd]^{\bundle{}fF{}}&\\
\Delta_N\ar[rr]^{\pi}&&N
}
\]

Theorem \ref{theorem: reduced fibered correspondence}
gives an opportunity to build two categories:
\begin{itemize}
\item \AddIndex{category of reduced fibered correspondences}
{category of reduced fibered correspondences};
its objects are
bundles over selected base and its morphisms are
reduced fibered correspondences
\item \AddIndex{category of fibered correspondences over diagonal}
{category of fibered correspondences over diagonal};
its objects are
bundles over selected base $M$ and its morphisms are
fibered correspondences base of which is diagonal in $M$
\end{itemize}
Functor between these categories is trivial. It maps
objects and morphisms into themselves, however in case of morphisms
we substitute a projection on base
by a projection on diagonal.

	\begin{remark}
	\label{remark: base of reduced fibered correspondences}
As evidenced by the foregoing, we can assume without loss of generality
that $N=M$. Like in
\ifx\PrintBook\Defined
remark \ref{remark: reduce details on the diagram, algebra bundle}
\else
section\xRef{0702.561}{section: Representation of Fibered F-Algebra}
\fi
we may not indicate the base on diagram of reduced fibered correspondences.
According to theorem \ref{theorem: reduced fibered correspondence}
it is not important for us whether we use the set
$M$ as base or we use the set $\Delta_M$.
		\qed
	\end{remark}

Let transition functions $f_{\alpha\beta}$ determine bundle $\mathcal{A}$
over base $M$ and transition functions $g_{\alpha\beta}$
determine bundle $\mathcal{B}$ over base $M$.
Let us consider
maps $U_\alpha\in M$ and $U_\beta\in M$, $U_\alpha\cap U_\beta\ne\emptyset$.
Point $p\in\mathcal{A}$ has representation $(x,p_\alpha)$
in map $U_\alpha$ and representation
$(x,p_\beta)$ in map $U_\beta$.
Point $q\in\mathcal{B}$ has representation $(x,q_\alpha)$
in map $U_\alpha$ and representation
$(x,q_\beta)$ in map $U_\beta$.
Therefore,
\[
p_\alpha=f_{\alpha\beta}(p_\beta)
\]
\[
q_\alpha=g_{\alpha\beta}(q_\beta)
\]
When we move from map $U_\alpha$ to map $U_\beta$,
representation of correspondence changes according to the law
\[
(x,p_\alpha,q_\alpha)=(x,f_{\alpha\beta}(p_\beta),g_{\alpha\beta}(q_\beta))
\]
This is consistent with the transformation
when we move from map $U_\alpha$ to map $U_\beta$
in the bundle $\mathcal{A}\odot\mathcal{B}$.

		\begin{theorem}
Given reduced
fibered correspondence $\mathcal{F}$ from $\mathcal{A}$ to $\mathcal{B}$,
set $\Gamma(\mathcal{F})$ determines correspondence
from $\Gamma(\mathcal{A})$ to $\Gamma(\mathcal{B})$.
		\end{theorem}
		\begin{proof}
By remark
\xRef{0702.561}{remark: reduced Cartesian product, section}
we can represent section of bundle $\mathcal{A}\odot\mathcal{B}$ as
tuple $(f,g)$ where $f$ is section of bundle $\mathcal{A}$
and $g$ is section of bundle $\mathcal{B}$. For each
fiber, $f(x)\in A_x$ is in
correspondence $F$ with $g(x)\in B_x$ iff section $(f,g)\in\Gamma(\mathcal{F})$.
		\end{proof}

The properties of reduced fibered correspondence
are closer to properties ordinary correspondence.

		\begin{theorem}
Let $\mathcal{F}$ be reduced fibered correspondence
from $\mathcal{A}$ to $\mathcal{B}$.
Suppose the bundle $\mathcal{A'}\subseteq\mathcal{A}$.
We define the image
of the bundle $\mathcal{A'}$ under reduced
fibered correspondence $\mathcal{F}$
according to law
\[
\mathcal{F}\mathcal{A}'=\{(x,b):x\in N,(a,b)\in F_x,a\in A_x\}
\]
Image of bundle $\mathcal{A'}$ under reduced fibered
correspondence $\mathcal{F}$
is subbundle of bundle $\mathcal{B}$.
		\end{theorem}
		\begin{proof}
Let us consider\footnote{I can repeat, up to notation,
proof of theorem
\ref{theorem: image of fibered correspondence}. However
I want to give another proof in order
to show how theorem
\ref{theorem: reduced fibered correspondence} works.}
commutative diagram of fibered correspondences
	\begin{equation}
\xymatrix{
\mathcal{A}'\ar[d]^{F/A'}\ar[rr]^I&&
\mathcal{A}\ar[d]^F\\
\mathcal{B}'\ar[rr]^J&&\mathcal{B}
}
	\label{eq: image of reduced fibered correspondence}
	\end{equation}

Depending on projections selected, this diagram represents either the
relationship between reduced fibered correspondences
over base $M$, or the
relationship between fibered correspondences over base $\Delta_M$.
However, these correspondences are identical.
All the relationships that hold for the base $\Delta_M$
hold for base $M$ as well.
		\end{proof}

	\begin{remark}
Diagram
\eqref{eq: image of reduced fibered correspondence} has simple representation.
However, if we draw this diagram without considering remark
\ref{remark: base of reduced fibered correspondences},
this diagram will has representation
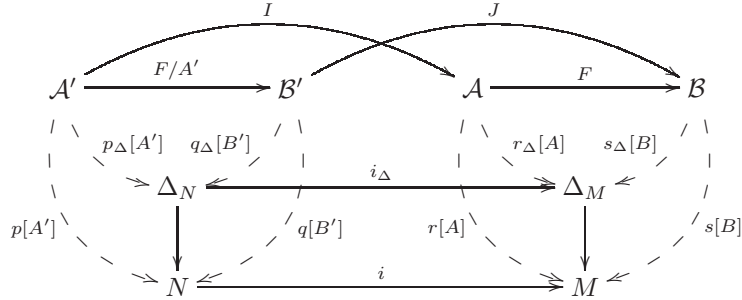
\[
\xymatrix{\\
\mathcal{A}'\ar[rr]^{F/A'}\ar@/^2pc/[rrrr]^I\arp@/_1pc/[dr]^{\bundle{}{p_\Delta}{A'}{}}
\arp@/_2pc/[ddr]_{\bundle{}p{A'}{}}&&
\mathcal{B}'\ar@/^2pc/[rrrr]^J\arp@/^1pc/[dl]_{\bundle{}{q_\Delta}{B'}{}}
\arp@/^2pc/[ddl]^{\bundle{}q{B'}{}}&&
\mathcal{A}\ar[rr]^F\arp@/_1pc/[dr]^{\bundle{}{r_\Delta}A{}}
\arp@/_2pc/[ddr]_{\bundle{}rA{}}&&
\mathcal{B}\arp@/^1pc/[dl]_{\bundle{}{s_\Delta}B{}}
\arp@/^2pc/[ddl]^{\bundle{}sB{}}\\
&\Delta_N\ar[d]\ar[rrrr]^{i_\Delta}&&&&\Delta_M\ar[d]\\
&N\ar[rrrr]^i&&&&M
}
\]
On diagram we use convention
\[
i_\Delta(a,a)=(i(a),i(a))
\]
		\qed
	\end{remark}

		\begin{theorem}
		\label{theorem: composition of reduced fibered correspondences}
Let us consider reduced fibered correspondence $\bundle{}sF{}$
from bundle $\mathcal{A}$ to bundle $\mathcal{B}$
and reduced fibered correspondence $\bundle{}tH{}$
from bundle $\mathcal{B}$ to bundle $\mathcal{C}$.
We define a
\AddIndex{composition of reduced fibered correspondences}
{composition of reduced fibered correspondences}
 $\bundle{}sF{}$ and $\bundle{}tH{}$
\[
\bundle{}tH{}\circ\bundle{}sF{}=\{(x,a,c):x\in M,(a,b)\in F_x,(b,c)\in H_x\}
\]
		\end{theorem}
		\begin{proof}
From commutativity of diagram of fibered correspondence
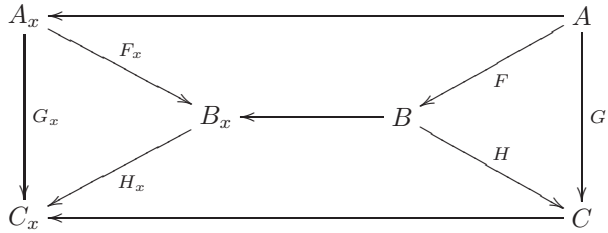
\[
\xymatrix{
A_x\ar[dd]^{G_x}\ar[drr]^{F_x}&&&&&&A\ar[dd]^G\ar[dll]^F\ar[llllll]\\
&&B_x\ar[dll]^{H_x}&&B\ar[drr]^H\ar[ll]&&\\
C_x&&&&&&C\ar[llllll]
}
\]
over base $\Delta_M$ follows $G_x=H_x\circ F_x$. Therefore
this equation holds also over base $M$.
		\end{proof}

		\begin{definition}
		\label{definition: inverse reduced fibered correspondence}
Let $\bundle{\mathcal{F}}sFM$ be reduced fibered correspondence
\[
\xymatrix{
\mathcal{A}\arp@/_1pc/[dr]_p\ar[rr]^{\mathcal{F}}&
\arp[d]^{\bundle{}sF{}}&\mathcal{B}\arp@/^1pc/[dl]^q\\
&M&
}
\]
from bundle $\mathcal{A}$ to bundle $\mathcal{B}$.
Then there exist
\AddIndex{inverse reduced fibered correspondence}
{inverse reduced fibered correspondence}
\symb{\mathcal{F}^{-1}}0{inverse reduced fibered correspondence, 1}
\symb{\bundle{}s{F^{-1}}{}}0{inverse reduced fibered correspondence, 2}
$\bundle{\ShowSymbol{inverse reduced fibered correspondence, 1}}s{F^{-1}}M$
\[
\xymatrix{
\mathcal{B}\arp@/_1pc/[drr]_p
\ar[rrrr]^{\ShowSymbol{inverse reduced fibered correspondence, 1}}&&
\arp[d]^{\ShowSymbol{inverse reduced fibered correspondence, 2}}&&
\mathcal{A}\arp@/^1pc/[dll]^q\\
&&M&&
}
\]
		\qed
		\end{definition}

		\begin{theorem}
Diagonal $\Delta_\mathcal{A}$ in bundle $\bundle{\mathcal{A}}pAM$
is reduced  fibered correspondence
\symb{\bundle{}r{\Delta_A}{}}0{diagonal in reduced bundle, 2}
$$\ShowSymbol{diagonal in reduced bundle, 2}:\Delta_\mathcal{A}\rightarrow M$$
where projection $\bundle{}r{\Delta_A}{}$ is determined by law
$$\bundle{}r{\Delta_A}{}(x,(p,p))=x$$
		\qed
		\end{theorem}

		\begin{theorem}
Let $\mathcal F$ be reduced fibered correspondence
from bundle $\mathcal{A}$ to bundle $\mathcal{B}$,
$\mathcal H$ be reduced fibered correspondence
from bundle $\mathcal{B}$ to bundle $\mathcal{C}$.
The following laws are valid for
reduced fibered correspondences $\mathcal F$ and
$\mathcal H$
\[
(\mathcal H\circ\mathcal F)^{-1}=\mathcal F^{-1}\circ\mathcal H^{-1}
\]
\[
(\mathcal F^{-1})^{-1}=\mathcal F
\]
\[
\mathcal F\circ\Delta_\mathcal{A}=\Delta_\mathcal{B}\circ\mathcal F=\mathcal F
\]
		\end{theorem}
		\begin{proof}
The proof of theorem is similar to proof of theorem
\ref{theorem: associative law, composition of fibered correspondences}.
We check each statement for base and fiber.
		\end{proof}

			\section{Fibered Relation}
			\label{section: Fibered Relation}

	\begin{definition}
	\label{definition: Fibered Relation}
Let $\bundle{\mathcal{A}}pAM$ be
bundle and $\omega$ be $n$\Hyph ary relation in the set $A$.
Fibered subset $\bundle{}r\omega{}$ of bundle
$\mathcal{E}^n$ is \AddIndex{$n$\Hyph ary fibered relation}
{fibered relation} in bundle $\mathcal{A}$.
	 \qed
	 \end{definition}

		\begin{theorem}
		\label{theorem: 2 ary fibered relation}
\AddIndex{$2$\Hyph ary fibered relation}
{2 ary fibered relation} in bundle $\mathcal{A}$ is
reduced fibered correspondence in bundle $\mathcal{A}$.
		\end{theorem}
		\begin{proof}
The theorem is corollary of definitions
\ref{definition: reduced fibered correspondence} and
\ref{definition: Fibered Relation}.
		\end{proof}

	\begin{definition}
$2$\Hyph ary fibered relation $\mathcal{F}$ in bundle $\mathcal{A}$
is said to be
\begin{itemize}
\item \AddIndex{transitive}
{transitive 2 ary fibered relation}, if
$\mathcal{F}\circ\mathcal{F}\subseteq\mathcal{F}$
\item \AddIndex{symmetric}{symmetric 2 ary fibered relation},
if $\mathcal{F}^{-1}=\mathcal{F}$
\item \AddIndex{antisymmetric}{antisymmetric 2 ary fibered relation},
if $\mathcal{F}\cap\mathcal{F}^{-1}\subseteq\Delta_\mathcal{A}$
\item \AddIndex{reflexive}{reflexive 2 ary fibered relation},
if $\mathcal{F}\supseteq\Delta_\mathcal{A}$
\end{itemize}
	 \qed
	 \end{definition}

	\begin{definition}
A transitive reflexive $2$\Hyph ary fibered relation $\mathcal{F}$
in bundle $\mathcal{A}$ is called a
\AddIndex{fibered preordering}
{fibered preordering} of $\mathcal{A}$.
$\mathcal{F}^{-1}$ is then also a fibered preordering
of $\mathcal{A}$;
it is said to be \AddIndex{opposite}
{opposite fibered preordering} to $\mathcal{F}$.
An antisymmetric fibered preordering
in bundle $\mathcal{A}$ is called a
\AddIndex{fibered ordering}
{fibered ordering} of $\mathcal{A}$.\footnote{One may be tempted to
define fibered total ordering $\mathcal{F}$
using equation $$\mathcal{F}\cup\mathcal{F}^{-1}=
\mathcal{A}^2$$ However, if we consider this relation on the set
of sections, then we can find two section which we cannot compare.}
	 \qed
	 \end{definition}

	\begin{definition}
A transitive reflexive symmetric
$2$\Hyph ary fibered relation $\mathcal{F}$
in bundle $\mathcal{A}$ is called a
\AddIndex{fibered equivalence}
{fibered equivalence} on bundle $\mathcal{A}$.
	 \qed
	 \end{definition}

\def\texFiberedMorphism{}
\ifx\PrintBook\Defined
				\chapter{Fibered Morphism}
\fi

			\section{Fibered Morphism}

		\begin{theorem}
		\label{theorem: quotient bundle}
Let us consider fibered equivalence $\bundle{\mathcal{S}}sSM$
on the bundle $\bundle{\mathcal{E}}pEM$.
Then there exists bundle
\[\bundle{\mathcal{E}/\mathcal{S}}t{E/S}M\]
called \AddIndex{quotient bundle}{quotient bundle}
of bundle $\mathcal{E}$ by
the equivalence $\mathcal{S}$.
Fibered morphism
\[
\mathrm{nat}\mathcal{S}:\mathcal{E}\rightarrow
\mathcal{E}/\mathcal{S}
\]
is called \AddIndex{fibered natural morphism}{fibered natural morphism}
or \AddIndex{fibered identification morphism}{fibered identification morphism}.
		\end{theorem}
		\begin{proof}
Let us consider the commutative diagram
\begin{equation}
\label{eq: quotient bundle}
\xymatrix{
\mathcal{E}\ar[rr]^{\mathrm{nat}\mathcal{S}}\arp[rd]_{\bundle{}pE{}}
&&\mathcal{E}/\mathcal{S}\arp[dl]^{\bundle{}t{E/S}{}}\\
&M&
}
\end{equation}
We introduce in $\mathcal{E}/\mathcal{S}$ quotient topology
(\citeBib{Bourbaki: General Topology 1}, page 33),
demanding continuity of mapping $\mathrm{nat}\mathcal{S}$.
According to proposition \citeBib{Bourbaki: General Topology 1}-I.3.6
mapping $\bundle{}t{E/S}{}$ is continuous.

Because we defined equivalence $S$ only between points of
the same fiber $E$, equivalence classes belong to
the same fiber $E/S$ (compare with the remark to proposition
\citeBib{Bourbaki: General Topology 1}-I.3.6).
		\end{proof}

Let $f:\mathcal{A}\rightarrow\mathcal{B}$ be
fibered morphism, base of which is
identity mapping.
According to definition
\xRef{0707.2246}{definition: inverse reduced fibered correspondence}
there exists inverse reduced fibered correspondence $f^{-1}$.
According to theorems
\xRef{0707.2246}{theorem: composition of reduced fibered correspondences}
and \xRef{0707.2246}{theorem: 2 ary fibered relation}
$f^{-1}\circ f$ is $2$\Hyph ary fibered relation.

		\begin{theorem}
Fibered relation $\mathcal{S}=f^{-1}\circ f$ is
fibered equivalence on the bundle $\mathcal{A}$.
There exists decomposition of fibered morphism $f$ into product of
fibered morphisms
	\begin{equation}
f=itj
	\label{eq: morphism of representations of algebra, homomorphism, 1}
	\end{equation}
\[
\xymatrix{
\mathcal{A}/\mathcal{S}\ar[rr]^t&&f(\mathcal{A})\ar[d]^i\\
\mathcal{A}\ar[u]^j\ar[rr]^f&&\mathcal{B}
}
\]
$j=\mathrm{nat}\mathcal{S}$ is the natural homomorphism
	\begin{equation}
\textcolor{red}{j(a)}=j(a)
	\label{eq: morphism of representations of algebra, homomorphism, 2}
	\end{equation}
$t$ is isomorphism
	\begin{equation}
\textcolor{red}{r(a)}=t(\textcolor{red}{j(a)})
	\label{eq: morphism of representations of algebra, homomorphism, 3}
	\end{equation}
$i $ is the inclusion mapping
	\begin{equation}
r(a)=i(\textcolor{red}{r(a)})
	\label{eq: morphism of representations of algebra, homomorphism, 4}
	\end{equation}
		\end{theorem}
		\begin{proof}
We verify the statement of theorem in fiber. We need also to check
that equivalence depends continuously on fiber.
		\end{proof}

\section{Free \texorpdfstring{$T\star$}{T*}-Representation of Fibered Group}

Mapping $\mathrm{nat}\mathcal{S}$ does not create bundle,
because different equivalence classes are not homeomorphic
in general.
However the proof of theorem \ref{theorem: quotient bundle}
suggests to the construction which reminds the construction
designed in \citeBib{Postnikov: Differential Geometry}, pages 16 - 17.

	 \begin{definition}
	 \label{definition: fibered little group}
Consider \Ts representation $f$ of fibered group $\bundle{}pG{}$
in bundle $\mathcal{M}$.
A \AddIndex{fibered little group}{fibered little group} or
\AddIndex{fibered stability group}{fibered stability group}
of $h\in \Gamma(\mathcal{M})$ is the set
\symb{\mathcal{G}_h}0{fibered little group}
\symb{\mathcal{G}_h}0{fibered stability group}
\[
\ShowSymbol{fibered little group}=\{g\in \Gamma(\mathcal{G}):f(g)h=h\}
\]
	 \qed
	 \end{definition}

	 \begin{definition}
	 \label{definition: free representation of fibered group}
\Ts representation $f$ of group $G$ is said to be
\AddIndex{free}{free representation of fibered group},
if for any $x\in M$ stability group $G_x=\{e\}$.
	 \qed
	 \end{definition}

		\begin{theorem}	%
		\label{theorem: free representation of group}
Given free \Ts representation $f$ of group $G$ in the set $A$,
there exist $1-1$ correspondence between orbits of representation,
as well between orbit of representation and group $G$.
		\end{theorem}
		\begin{proof}
		\end{proof}

Let us consider covariant free \Ts representation $f$
of fibered group $\bundle{}pG{}$ in fiber $\bundle{}pE{}$.
This \Ts representation determines fibered
equivalence $\mathcal{S}$ on $\bundle{}aE{}$,
$(p,q)\in S$ when $p$ and $q$ belong to
common orbit.
Since the representation in every fiber is free,
all equivalence classes are homeomorphic to group $G$.
Therefore, the mapping $\mathrm{nat}\mathcal{S}$
is projection of the bundle
$\bundle{\mathcal{E}}{\mathrm{nat}\mathcal{S}}G{\mathcal{E}/\mathcal{S}}$.
We also use notation $\mathcal{S}=\mathcal{G}\star$.
We may represent diagram \eqref{eq: quotient bundle}
in the following form
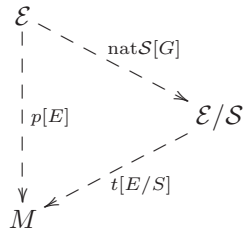
\[
\xymatrix{
\mathcal{E}\arp[drr]^{\bundle{}{\mathrm{nat}\mathcal{S}}G{}}
\arp[dd]^{\bundle{}pE{}}&&\\
&&\mathcal{E}/\mathcal{S}\arp[dll]^{\bundle{}t{E/S}{}}\\
M&&
}
\]
Bundle $\bundle{}{\mathrm{nat}\mathcal{S}}G{}$ is called
\AddIndex{bundle of level $2$}{bundle of level 2}.

\begin{example}
Let us consider the representation of rotation group $SO(2)$ in $R^2$.
All points except the point
$(0,0)$ have trivial little group.
Hence, we defined free representation of group $SO(2)$
in set $R^2\setminus \{(0,0)\}$.

We cannot use this idea in case of bundle
$\bundle{}p{R^2}{}$ and representation of fibered group $\bundle{}t{SO(2)}{}$.
Let $\mathcal{S}$ be relation of fibered equivalence.
The bundle $\bundle{}p{R^2\setminus \{(0,0)\}}{}/\bundle{}t{SO(2)}{}\star$
is not complete. As a consequence passage to the limit may bring into
non\Hyph existent fiber. Therefore we prefer to consider
bundle $\bundle{}p{R^2}{}/\bundle{}t{SO(2)}{}\star$,
keeping in mind, that fiber over point $(x,0,0)$ is degenerate.
	 \qed
\end{example}

We simplify the notation and represent this construction as
\symb{p[E_2,E_1]}0{bundle of level 2}
\[
\xymatrix{
\ShowSymbol{bundle of level 2}:\mathcal{E}_2\arp[r]&
\mathcal{E}_1\arp[r]&
M
}
\]
where we consider bundles
\begin{align*}
\bundle{\mathcal{E}_2}{p_2}{E_2}{\mathcal{E}_1}&&
\bundle{\mathcal{E}_1}{p_1}{E_1}M&
\end{align*}
Similarly we consider
\AddIndex{bundle of level $n$}{bundle of level n}
\symb{p[E_n,...,E_1]}0{bundle of level n}
\begin{equation}
\label{eq: bundle of level n}
\xymatrix{
\ShowSymbol{bundle of level n}:\mathcal{E}_n\arp[r]&\ar@{}[r]_{...}&\arp[r]&
\mathcal{E}_1\arp[r]&
M
}
\end{equation}
The sequence of bundles \eqref{eq: bundle of level n}
is called \AddIndex{tower of bundles}{tower of bundles}.
I made this definition by analogy with Postnikov tower
(\citeBib{Hatcher: Algebraic Topology}).
Postnikov tower is the tower of bundles. Fiber of bundle of level $n$ is
homotopy group of oreder $n$.
Such definitions are well known, however I gave definition of tower of
bundles, because it follows in a natural way from the text above.

One more example of tower of bundles
attracted my attention (\citeBib{geometry of differential equations},
\citeBib{cohomological analysis}, chapter 2).
We consider the set $J^0(n,m)$ of $0$\Hyph jets
of functions from $R^n$ to $R^m$ as base. We consider 
the set $J^p(n,m)$ of $p$\Hyph jets
of functions from $R^n$ to $R^m$ as bundle of level $p$.

\OpenBiblio


\BiblioItem{texIntro}{Einstein: Geometry and Experience}
{
Einstein, Geometry and Experience, (1921)
}%

\BiblioItem{texGenRelativity}{Ghez}
{
Ghez et al.,
The First Measurement of Spectral Lines in a Short-Period Star Bound to the Galaxy's Central Black Hole: A Paradox of Youth,
\href{http://www.journals.uchicago.edu/ApJ/journal/issues/ApJL/v586n2/16990/brief/16990.abstract.html}{ApJL, 586, L127} (2003),
eprint \href{http://arxiv.org/abs/astro-ph/0302299}{arXiv:astro-ph/0302299} (2003)
}%

\BiblioItem{texGenRelativity}{Schodel}
{
R. Sch\"odel et al.,
A star in a 15.2-year orbit around the supermassive black hole at the centre of the Milky Way,
\href{http://www.nature.com/cgi-taf/DynaPage.taf?file=/nature/journal/v419/n6908/abs/nature01121_fs.html}{Nature 419, 694} (2002)
}%

\BiblioItem{texAffine}{Mielke}
{
Eckehard W. Mielke, Affine generalization of the Komar complex of general relativity,
\href{http://prola.aps.org/searchabstract/PRD/v63/i4/e044018}{Phys. Rev. D 63, 044018} (2001)
}%

\BiblioItem{texAffine}{Obukhov}
{
Yu. N. Obukhov and J. G. Pereira, Metric-affine approach to teleparallel gravity,
\href{http://scitation.aip.org/getabs/servlet/GetabsServlet?prog=normal&id=PRVDAQ000067000004044016000001&idtype=cvips&gifs=Yes}
{Phys. Rev. D 67, 044016} (2003),
eprint \href{http://arxiv.org/abs/gr-qc/0212080}{arXiv:gr-qc/0212080} (2002)
}%

\BiblioItem{texAffine}{Sardanashvily}
{
Giovanni Giachetta, Gennadi Sardanashvily, Dirac Equation in Gauge and Affine-Metric Gravitation Theories,
eprint \href{http://arxiv.org/abs/gr-qc/9511035}{arXiv:gr-qc/9511035} (1995)
}%

\BiblioItem{texAffine}{Gauge}
{
Frank Gronwald and Friedrich W. Hehl, On the Gauge Aspects of Gravity, eprint
\href{http://arxiv.org/abs/gr-qc/9602013}{arXiv:gr-qc/9602013} (1996)
}%

\BiblioItem{texAffine}{Neeman}
{
Yuval Neeman, Friedrich W. Hehl, Test Matter in a Spacetime with Nonmetricity, eprint
\href{http://arxiv.org/abs/gr-qc/9604047}{arXiv:gr-qc/9604047} (1996)
}%

\BiblioItem{texAffine}{0405.027}
{
Aleks Kleyn,
Reference Frame in General Relativity,
eprint \href{http://arxiv.org/abs/gr-qc/0405027}{arXiv:gr-qc/0405027} (2004)
}%

\BiblioItem{texTidal,texAffine}{torsion}
{
F. W. Hehl, P. von der Heyde, G. D. Kerlick, and J. M. Nester,
General relativity with spin and torsion: Foundations and prospects,
\href{http://prola.aps.org/abstract/RMP/v48/i3/p393_1}{Rev. Mod. Phys. 48, 393} (1976)
}%

\BiblioItem{texTidal,texNewton}{Megged}
{
O. Megged, Post-Riemannian Merger of Yang-Mills Interactions with Gravity,
eprint \href{http://arxiv.org/abs/hep-th/0008135}{arXiv:hep-th/0008135} (2001)
}%


\BiblioItem{texNewton}{gr-qc-9604027}
{
Yu.N. Obukhov, E.J. Vlachynsky, W. Esser, R. Tresguerres and F.W. Hehl,
An exact solution of the metric-affine gauge theory with dilation, shear, and spin charges,
eprint \href{http://arxiv.org/abs/gr-qc/9604027}{arXiv:gr-qc/9604027} (1996)
}%

\BiblioItem{texLagrange}{Weinberg}
{
Steven Weinberg. The Quantum Theory of Fields. Cambridge university press.
}%

\BiblioItem{texLagrange}{Reinhardt}
{
Greiner Reinhardt. Field Quantization. Springer.
}%

\BiblioItem{texLagrange}{Landau}
{
L. D. Landau, E. M. Lifshich, The classical theory of fields.
Oxford, New York, Pergamon Press
}%

\BiblioItem{texTidal}{Wheeler}
{
Ignazio Ciufolini, John Wheeler. Gravitation and Inertia.
Princeton university press.
}%

\BiblioItem{texTidal}{0405.028}
{
Aleks Kleyn, Metric-Affine Manifold,
eprint \href{http://arxiv.org/abs/gr-qc/0405028}{arXiv:gr-qc/0405028} (2004)
}%

\BiblioItem{texTidal}{Anderson02}
{
J. D. Anderson, P. A. Laing, E. L. Lau, A. S. Liu, M. M. Nieto, and S. G. Turyshev,
Study of the anomalous acceleration of Pioneer 10 and 11,
\href{http://prola.aps.org/searchabstract/PRD/v65/i8/e082004}{Phys. Rev. D 65, 082004, 50 pp.}, (2002),
eprint \href{http://arxiv.org/abs/gr-qc/0104064}{arXiv:gr-qc/0104064} (2001)
}%

\BiblioItem{texTidal}{Anderson98}
{
J. D. Anderson, P. A. Laing, E. L. Lau, A. S. Liu, M. M. Nieto, and S. G. Turyshev,
Indication, from Pioneer 10/11, Galileo, and Ulysses Data, of an Apparent Anomalous, Weak, Long-Range Acceleration,
\href{http://prola.aps.org/abstract/PRL/v81/i14/p2858_1}{Phys. Rev. Lett. 81, 2858}, (1998),
eprint \href{http://arxiv.org/abs/gr-qc/9808081}{arXiv:gr-qc/9808081} (1998)
}%


\BiblioItem{texReferenceFrame,texFiberedAlgebra}{Serge Lang}
{
Serge Lang,
Algebra, Springer, 2002
}%

\BiblioItem{texFiberedAlgebra,texTstarMorphism}{Burris Sankappanavar}
{
S. Burris, H.P. Sankappanavar,
A Course in Universal Algebra, Springer-Verlag (March, 1982),
\\eprint
\href{http://www.math.uwaterloo.ca/~snburris/htdocs/ualg.html}
{http://www.math.uwaterloo.ca/~snburris/htdocs/ualg.html}
\\(The Millennium Edition)
}%


\BiblioItem{texAffine,texRepresentation,texBasis,texDrcBasis,texLinearMap,texVectorSpace}{Rashevsky}
{
P. K. Rashevsky, Riemann Geometry and Tensor Calculus,\\
Moscow, Nauka, 1967
}%

\BiblioItem{texDrcBasis,texBasis}{Korn}
{
Granino A. Korn, Theresa M. Korn,
Mathematical Handbook for Scientists and Engineer,
McGraw-Hill Book Company, New York, San Francisco,
Toronto, London, Sydney, 1968
}%


\BiblioItem{texGenRelativity}{Tartaglia}
{
Angelo Tartaglia and Matteo Luca Ruggiero,
Angular Momentum Effects in Michelson\Hyph Morley Type Experiments,
Gen.Rel.Grav. 34, 1371-1382 (2002),\\
eprint \href{http://arxiv.org/abs/gr-qc/0110015}{arXiv:gr-qc/0110015} (2001)
}%

\BiblioItem{texGenRelativity}{Tomozawa}
{
Yukio Tomozawa, Speed of Light in Gravitational Fields, eprint
\href{http://arxiv.org/abs/astro-ph/0303047}{arXiv:astro-ph/0303047} (2004)
}%

\BiblioItem{texGenRelativity}{Magueijo}
{
Joao Magueijo,
Covariant and locally Lorentz-invariant varying speed of light theories,
\href{http://prola.aps.org/abstract/PRD/v62/i10/e103521}{Phys. Rev. D 62, 103521} (2000),
eprint \href{http://arxiv.org/abs/gr-qc/0007036}{arXiv:gr-qc/0007036} (2000)
}%

\BiblioItem{texGenRelativity}{Bassett}
{
Bruce A. Bassett, Stefano Liberati, Carmen Molina-Paris, and Matt Visser,
Geometrodynamics of variable-speed-of-light cosmologies,
\href{http://prola.aps.org/abstract/PRD/v62/i10/e103518}{Phys. Rev. D 62}, 103518 (2000),
eprint \href{http://arxiv.org/abs/astro-ph/0001441}{arXiv:astro-ph/0001441} (2000)
}%

\BiblioItem{texGenRelativity}{Straumann}
{
Lochlainn O'Raifeartaigh and Norbert Straumann,
Gauge theory: Historical origins and some modern developments,
\href{http://prola.aps.org/abstract/RMP/v72/i1/p1_1}{Rev. Mod. Phys. 72, 1} (2000)
}%

\BiblioItem{texGenRelativity}{Lammerzahl}
{
Claus L\"ammerzahl, Mark P. Haugan,
On the interpretation of Michelson\Hyph Morley experiments,
{Phys. Lett. A282 223-229} (2001),\\
eprint \href{http://arxiv.org/abs/gr-qc/0103052}{arXiv:gr-qc/0103052} (2001)
}%

\BiblioItem{texGenRelativity}{Muller}
{
Holger Muller et al.,
Modern Michelson-Morley Experiment using Cryogenic Optical Resonators,
\href{http://prola.aps.org/searchabstract/PRL/v91/i2/e020401}{Phys. Rev. Lett. 91, 020401} (2003),
eprint \href{http://arxiv.org/abs/physics/0305117}{arXiv:physics/0305117} (2000)
}%

\BiblioItem{texGenRelativity,texTidal}{Ranada}
{
Antonio F. Ranada,
Pioneer acceleration and variation of light speed: experimental situation,
eprint \href{http://arxiv.org/abs/gr-qc/0402120}{arXiv:gr-qc/0402120} (2004)
}%

\BiblioItem{texBiring,texVectorSpace}{math.QA-0208146}
{
I. Gelfand, S. Gelfand, V. Retakh, R. Wilson,
Quasideterminants,\\
eprint \href{http://arxiv.org/abs/math.QA/0208146}{arXiv:math.QA/0208146} (2002)
}%

\BiblioItem{texBiring,texVectorSpace}{q-alg-9705026}
{
I.Gelfand, V.Retakh,
Quasideterminants, I,\\
eprint \href{http://arxiv.org/abs/q-alg/9705026}{arXiv:q-alg/9705026} (1997)
}%

\BiblioItem{texVectorSpace}{Gelfand Retakh 1991}
{
I. Gelfand and V. Retakh, Determinants of Matrices over Noncommutative Rings, Funct.
Anal. Appl. 25 (1991), no. 2, 91-102
}%

\BiblioItem{texVectorSpace}{Gelfand Retakh 1992}
{
I. Gelfand and V. Retakh, A Theory of Noncommutative Determinants and Characteristic
Functions of Graphs, Funct. Anal. Appl. 26 (1992), no. 4, 1-20
}%

\BiblioItem{texVectorSpace}{hep-th-9407124}
{
I. M. Gelfand, D. Krob, A. Lascoux, B. Leclerc, V.S. Retakh and J.-Y. Thibon,
Noncommutative symmetric functions,\\
eprint \href{http://arxiv.org/abs/hep-th/9407124}{arXiv:hep-th/9407124} (1994)
}%

\BiblioItem{texVectorSpace}{Carl Faith 1}
{
Carl Faith, Algebra: Rings, Modules and Categories I,
Springer - Verlag, Berlin - Heidelberg - New York, 1973
}%



\BiblioItem{texReferenceFrame}{math.DG-0412391}
{
Aleks Kleyn,
Basis Manifold,
eprint \href{http://arxiv.org/abs/math.DG/0412391}{arXiv:math.DG/0412391} (2004)
}%

\BiblioItem{texFiberedAlgebra,texBundleRelation,texTstarMorphism}{0701.238}
{
Aleks Kleyn,
Lectures on Linear Algebra over Skew Field,\\
eprint \href{http://arxiv.org/abs/math.GM/0701238}{arXiv:math.GM/0701238} (2007)
}%

\BiblioItem{texBundleRelation,texPrefaceRelation}{0702.561}
{
Aleks Kleyn,
Algebra Bundle,\\
eprint \href{http://arxiv.org/abs/math.DG/0702561}{arXiv:math.DG/0702561} (2007)
}%


\BiblioItem{texPolymodule}{math.RA-0501237v1}
{
Aleks Kleyn,
Module Over Skew-Field, version 1,\\
eprint \href{http://arxiv.org/abs/math/0501237v1}{arXiv:math.RA/0501237v1} (2005)
}%

\ifx\texBiring\Defined
\else
\BiblioItem{texVectorSpace,texFiberedAlgebra}{0612.111}
{
Aleks Kleyn,
Biring of Matrices,\\
eprint \href{http://arxiv.org/abs/math.OA/0612111}{arXiv:math.OA/0612111} (2006)
}%
\fi

\ifx\texBundleRelation\Defined
\else
\BiblioItem{texFiberedMorphism}{0707.2246}
{
Aleks Kleyn,
Fibered Correspondence,\\
eprint \href{http://arxiv.org/abs/0707.2246}{arXiv:0707.2246} (2007)
}%
\fi

\BiblioItem{texHomotopy}{q-alg-9705009}
{
John C. Baez,
An Introduction to n-Categories,\\
eprint \href{http://arxiv.org/abs/q-alg/9705009}{arXiv:q-alg/9705009} (1997)
}%

\BiblioItem{texIntro}{Einstein: Isaak Newton}
{
Albert Einstein,
Isaak Newton,
Manchester Guardian, 16, 234 - 235 (1927)
}%

\BiblioItem{texPrefaceRelation}{Tolstoi about Anna Karenina}
{
Tolstoi about Anna Karenina,
in book A Karenina Companion, by C. J. G. Turner,
published by Wilfrid Laurier University Press (August 1993)
}%

\BiblioItem{texBundleRelation,texPrefaceRelation,texTstarMorphism,texCartesian}
{Cohn: Universal Algebra}
{
Paul M. Cohn,
Universal Algebra,
Springer, 1981
}%

\BiblioItem{texCartesian}
{Maunder: Algebraic Topology}
{
C. R. F. Maunder,
Algebraic Topology,
Dover Publications, Inc, Mineola, New York, 1996
}%

\BiblioItem{texFiberedAlgebra}{Pommaret: Partial Differential Equations}
{
J.-F. Pommaret,
Partial Differential Equations and Group Theory,
Springer, 1994
}%

\BiblioItem{texBundleRelation}{Bourbaki: Set Theory}
{
N. Bourbaki,
Theory of sets,
Springer, 2004
}%

\BiblioItem{texCartesian,texFiberedAlgebra,texBundleRelation,texFiberedMorphism}
{Bourbaki: General Topology 1}
{
N. Bourbaki,
General Topology, Chapters 1 - 4,
Springer, 1989
}

\BiblioItem{texCalculus}{Bourbaki: General Topology: Chapter 5 - 10}
{
N. Bourbaki,
General Topology, Chapters 5 - 10,
Springer, 1989
}

\BiblioItem{texCalculus}{Bourbaki: Topological Vector Space}
{
N. Bourbaki,
Topological Vector Spaces, Chapters 1 - 5,
Springer, 1987
}

\BiblioItem{texCalculus}{Pontryagin: Topological Group}
{
L. S. Pontryagin,
Selected Works, Volume Two, Topological Groups,
Gordon and Breach Science Publishers, 1986
}

\BiblioItem{texFiberedMorphism}{Postnikov: Differential Geometry}
{
Postnikov M. M.,
Geometry IV: Differential geometry,
Moscow, Nauka, 1983
}

\BiblioItem{texFiberedAlgebra,texFiberedMorphism}{Hatcher: Algebraic Topology}
{
Allen Hatcher,
Algebraic Topology,
Cambridge University Press, 2002
}

\BiblioItem{texFiberedMorphism}{geometry of differential equations}
{
Vinogradov, A. M., Krasil'shchik, I. S., and Lychagin, V. V.,
Introduction to geometry of nonlinear differential equations,
Nauka, Moscow, 1986
}

\BiblioItem{texFiberedMorphism}{cohomological analysis}
{
A. M. Vinogradov,
Cohomological Analysis of Partial Differential Equations
and Secondary Calculus,
American Mathematical Society, 2001
}

\CloseBiblio

\OpenIndex
\SetIndexSpace%
\Index{texLinearMap}
   {$1$-\drc form}%
   {1-drc form, vector spaces}%
\SetIndexSpace%
\Index{texPolymodule}
   {$(2)$\hyph vector space}%
   {(2)-vector space}%
\Index{texBundleRelation}
   {$2$\Hyph ary fibered relation}%
   {2 ary fibered relation}%
\SetIndexSpace%
\Index{texCalculus}
   {$A$\Hyph valued function}%
   {A valued function}%
\Index{texBiring}
   {$(^{\gi a}_{\gi b})$\hyph \CR quasideterminant}%
   {a b cr-quasideterminant}%
\Index{texBiring}
   {$(^{\gi a}_{\gi b})$\hyph \RC quasideterminant}%
   {a b RC-quasideterminant}%
\Index{texCalculus}
   {absolute value on skew field}%
   {absolute value on skew field}%
\Index{texBasis}
   {active representation}%
   {active representation}%
\Index{texDrcBasis}
   {active \sT representation}%
   {active representation, vector space}%
\Index{texBasis}
   {active transformation on basis manifold}%
   {active transformation}%
\Index{texBasis}
   {affine basis}%
   {Affine Basis}%
\Index{texBasis}
   {affine transformation group}%
   {AffineTransformationGroup}%
\Index{texTypeBasis}
   {affine transformation group}%
   {AffineTransformationGroup}%
\Index{texBasis}
   {affine transformation on basis manifold}%
   {affine transformation}%
\Index{texBiring}
   {alternative representation of matrix}%
   {Alternative representation}%
\Index{texReferenceFrame}
   {anholonomic coordinate}%
   {anholonomic coordinate}%
\Index{texReferenceFrame}
   {anholonomic coordinates of vector}%
   {vector anholonomic coordinates}%
\Index{texReferenceFrame}
   {anholonomic coordinates on manifold}%
   {anholonomic coordinates on manifold}%
\Index{texReferenceFrame}
   {anholonomity object}%
   {anholonomity object}%
\Index{texFiberedGroup}
   {antihomomorphism of fibered groups}%
   {antihomomorphism of fibered groups}%
\Index{texBundleRelation}
   {antisymmetric $2$\Hyph ary fibered relation}%
   {antisymmetric 2 ary fibered relation}%
\Index{texFiberedAlgebra}
   {arity of operation}%
   {arity of operation}%
\Index{texTstarRepresentation}
   {associative law for covariant \Ts representation}%
   {associative law for Tstar covariant representation}%
\Index{texDrcMorphism}
   {associative law for \drc linear maps of vector spaces}%
   {associative law for drc linear maps of vector spaces}%
\Index{texVectorSpace}
   {associative law for \Ts vector space}%
   {associative law, Tstar vector space}%
\Index{texLinearMap}
   {associative law for twin representations}%
   {associative law for twin representations}%
\Index{texBundleRelation}
   {associative law of composition of fibered correspondences}%
   {associative law, composition of fibered correspondences}%
\Index{texAffine}
   {auto parallel line}%
   {auto parallel line}%
\SetIndexSpace%
\Index{texBundleRelation}
   {base of fibered correspondence}%
   {base of fibered correspondence}%
\Index{texCartesian}
   {base of map}%
   {base of map}%
\Index{texTypeBasis}
   {basis}%
   {}%
\SubIndex{texTypeBasis}
   {affine}%
   {Affine Basis}%
\SubIndex{texTypeBasis}
   {central affine}%
   {Central Affine Basis}%
\SubIndex{texTypeBasis}
   {orthonornal}%
   {Orthonornal Basis}%
\Index{texDrcBasis}
   {basis manifold}%
   {}%
\SubIndex{texTypeBasis}
   {of affine space}%
   {Basis Manifold, Affine Space}%
\SubIndex{texTypeBasis}
   {of central affine space}%
   {Basis Manifold, Central Affine Space}%
\SubIndex{texDrcBasis}
   {of \drc vector space}%
   {basis manifold of vector space}%
\SubIndex{texTypeBasis}
   {of Euclid space}%
   {Basis Manifold, Euclid Space}%
\Index{texBasis}
   {basis manifold of affine space}%
   {Basis Manifold, Affine Space}%
\Index{texBasis}
   {basis manifold of central affine space}%
   {Basis Manifold, Central Affine Space}%
\Index{texBasis}
   {basis manifold of Euclid space}%
   {Basis Manifold, Euclid Space}%
\Index{texBasis}
   {basis manifold of vector space}%
   {basis manifold of vector space}%
\Index{texBasis}
   {basis of vector space}%
   {Basis}%
\Index{texLieRepresentation}
   {basis vector}%
   {}%
\SubIndex{texLieRepresentation}
   {of \sT representation}%
   {basis vector of starT representation}%
\SubIndex{texLieRepresentation}
   {of \Ts representation}%
   {basis vector of Tstar representation}%
\Index{texBiring}
   {biring}%
   {biring}%
\Index{texFiberedMorphism}
   {bundle of level $2$}%
   {bundle of level 2}%
\Index{texFiberedMorphism}
   {bundle of level $n$}%
   {bundle of level n}%
\SetIndexSpace%
\Index{texVectorSpace}
   {\subs rows \dcr vector space}%
   {subs rows dcr vector space}%
\Index{texBiring}
   {\subs row of matrix}%
   {c row}%
\Index{texBiring}
   {$c$\hyph row of matrix}%
   {c-row}%
\Index{texAffine}
   {Cartan connection}%
   {Cartan connection}%
\Index{texAffine}
   {Cartan curvature}%
   {Cartan curvature}%
\Index{texAffine}
   {Cartan derivative}%
   {Cartan derivative}%
\Index{texAffine}
   {Cartan symbol}%
   {Cartan symbol}%
\Index{texAffine}
   {Cartan transport}%
   {Cartan transport}%
\Index{texCartesian}
   {Cartesian power $A$ of set $B$}%
   {Cartesian power of set}%
\Index{texCartesian}
   {Cartesian power of bundle}%
   {Cartesian power of bundle}%
\Index{texCartesian}
   {direct product of bundles}%
   {Cartesian product of bundles}%
\Index{texCartesian}
   {direct product of total spaces}%
   {Cartesian product of total spaces}%
\Index{texHomotopy}
   {category of \drc vector spaces}%
   {category of drc vector spaces}%
\Index{texBundleRelation}
   {category of fibered correspondences over diagonal}%
   {category of fibered correspondences over diagonal}%
\Index{texBundleRelation}
   {category of reduced fibered correspondences}%
   {category of reduced fibered correspondences}%
\Index{texBasis}
   {central affine basis}%
   {Central Affine Basis}%
\Index{texRepresentation}
   {column vector}%
   {column vector}%
\Index{texBundleRelation}
   {commutative diagram of correspondences}%
   {commutative diagram of correspondences}%
\Index{texCartesian}
   {compact\hyph open topology}%
   {compact open topology}%
\Index{texDiffEq}
   {completely integrable system}%
   {completely integrable system}%
\Index{texBundleRelation}
   {composition of fibered correspondences}%
   {composition of fibered correspondences}%
\SubIndex{texVectorSpace}
   {of \drc linear equations}%
   {extended matrix, system of drc linear equations}%
\SubIndex{texVectorSpace}
   {of \rcd linear equations}%
   {extended matrix, system of rcd linear equations}%
\Index{texBundleRelation}
   {composition of reduced fibered correspondences}%
   {composition of reduced fibered correspondences}%
\Index{texBiring}
   {condition of reducibility of products}%
   {condition of reducibility of products}%
\Index{texBundleRelation}
   {continuous correspondence}%
   {continuous correspondence}%
\Index{texFiberedGroup}
   {contravariant \Ts representation of fibered group}%
   {Tstar contravariant representation of fibered group}%
\Index{texBasis}
   {coordinate isomorphism}%
   {coordinate isomorphism}%
\Index{texVectorSpace}
   {coordinate isomorphism}%
   {coordinate isomorphism}%
\Index{texVectorSpace}
   {coordinate matrix}%
   {}%
\SubIndex{texVectorSpace}
   {of set of vectors in \dcr rows vector space}%
   {coordinate matrix of set of vectors, dcr vector space}%
\SubIndex{texVectorSpace}
   {of set of vectors in \drc rows vector space}%
   {coordinate matrix of set of vectors, drc vector space}%
\SubIndex{texVectorSpace}
   {of vector in \drc basis}%
   {coordinate matrix of vector in drc basis}%
\Index{texReferenceFrame}
   {coordinate reference frame}%
   {coordinate reference frame}%
\Index{texGenRelativity}
   {coordinate reference frame in event space}%
   {coordinate reference frame in event space}%
\Index{texDrcBasis}
   {coordinate representation in \drc vector space}%
   {coordinate representation, vector space}%
\Index{texBasis}
   {coordinate representation of group in vector space}%
   {coordinate representation, vector space}%
\Index{texBasis}
   {coordinate vector space}%
   {coordinate vector space}%
\Index{texDrcBasis}
   {coordinates of geometrical object}%
   {}%
\SubIndex{texDrcBasis}
   {in coordinate vector space}%
   {coordinates of geometrical object, coordinate vector space}%
\SubIndex{texDrcBasis}
   {in vector space}%
   {coordinates of geometrical object, vector space}%
\Index{texBasis}
   {coordinates of geometrical object}%
   {coordinates of geometrical object, vector space}%
\Index{texBasis}
   {coordinates of geometrical object in coordinate representation}%
   {coordinates of geometrical object, coordinate vector space}%
\Index{texDrcBasis}
   {coordinates of representation}%
   {coordinates of representation}%
\Index{texBasis}
   {coordinates of representation}%
   {coordinates of representation}%
\Index{texVectorSpace}
   {coordinates of set of vectors in \dcr vector space}%
   {coordinates of set of vectors, dcr vector space}%
\Index{texVectorSpace}
   {coordinates of set of vectors in \drc vector space}%
   {coordinates of set of vectors, drc vector space}%
\Index{texVectorSpace}
   {coordinates of vector in \drc basis}%
   {coordinates of vector in drc basis}%
\Index{texBundleRelation}
   {correspondence continuous on the set}%
   {correspondence continuous on the set}%
\Index{texBundleRelation}
   {correspondence of homomorphism}%
   {correspondence of homomorphism}%
\Index{texFiberedGroup}
   {covariant \Ts representation of fibered group}%
   {Tstar covariant representation of fibered group}%
\Index{texBiring}
   {\CR inverse element of biring}%
   {cr-inverse element}%
\Index{texVectorSpace}
   {\CR matrix group}%
   {cr-matrix group}%
\Index{texBiring}
   {\CR power}%
   {cr power}%
\Index{texBiring}
   {\CR product of matrices}%
   {cr-product of matrices}%
\Index{texVectorSpace}
   {\crd vector space}%
   {crd vector space}%
\SetIndexSpace%
\Index{texCalculus}
   {$D$\Hyph valued variable}%
   {D valued variable}%
\Index{texVectorSpace}
   {\dcr basis of \subs rows vector space}%
   {dcr basis, c rows vector space}%
\Index{texVectorSpace}
   {\dcr vector}%
   {dcr vector}%
\Index{texVectorSpace}
   {\dcr vector space}%
   {dcr vector space}%
\Index{texBiring}
   {determinant of matrix}%
   {determinant}%
\Index{texTidal}
   {deviation of trajectories}%
   {deviation of trajectories}%
\Index{texBundleRelation}
   {diagonal in bundle}%
   {diagonal in bundle}%
\Index{texBundleRelation}
   {diagram of correspondences}%
   {diagram of correspondences}%
\Index{texCalculus}
   {differentiable functions of \drc vector space to skew field $D$}%
   {differentiable functions, drc vector space to skew field}%
\Index{texCalculus}
   {differential of mapping of normed \drc vector space to valued skew field}%
   {differential, drc vector space to skew field}%
\Index{texVectorSpace}
   {dimension of \drc vector space}%
   {dimension of vector space}%
\Index{texFiberedGroup}
   {direct product of representations of fibered group}%
   {direct product of representations of fibered group}%
\Index{texRepresentation}
   {direct product of representations of group}%
   {direct product of representations of group}%
\Index{texTstarRepresentation}
   {direct product of \Ts representations of group}%
   {direct product of representations of group}%
\Index{texLieRepresentation}
   {direct sum of representations}%
   {direct sum of representations}%
\Index{texVectorSpace}
   {distributive law}%
   {}%
\SubIndex{texVectorSpace}
   {\Ts vector space}%
   {distributive law, Tstar vector space}%
\Index{texVectorSpace}
   {\drc automorphism of vector space}%
   {automorphism of vector space}%
\Index{texVectorSpace}
   {\drc basis}%
   {}%
\SubIndex{texVectorSpace}
   {for \sups rows vector space}%
   {drc basis, r rows vector space}%
\SubIndex{texVectorSpace}
   {for vector space}%
   {drc basis, vector space}%
\Index{texVectorSpace}
   {\drc coordinate vector space}%
   {drc coordinate vector space}%
\Index{texVectorSpace}
   {\drc isomorphism of vector spaces}%
   {isomorphism of vector spaces}%
\Index{texDrcMorphism}
   {\drc linear map of vector spaces}%
   {drc linear map of vector spaces}%
\Index{texVectorSpace}
   {\drc linear span in vector space}%
   {linear span, vector space}%
\Index{texVectorSpace}
   {\drc linearly dependent vectors}%
   {linearly dependent, vector space}%
\Index{texReferenceFrame}
   {\drc linearly independent vector fields}%
   {linearly independent vector fields}%
\Index{texVectorSpace}
   {\drc linearly independent vectors}%
   {linearly independent, vector space}%
\Index{texVectorSpace}
   {\drc system of linear equations}%
   {system of linear equations}%
\Index{texVectorSpace}
   {\drc vector}%
   {drc vector}%
\Index{texCalculus}
   {\drc vector function}%
   {drc vector function}%
\Index{texVectorSpace}
   {\drc vector space}%
   {drc vector space}%
\Index{texVectorSpace}
   {$D\star$\hyph  vector space}%
   {Dstar vector space}%
\Index{texVectorSpace}
   {$D\star$\hyph product of vector over scalar}%
   {Dstar product of vector over scalar, vector space}%
\Index{texBiring}
   {duality principle for biring}%
   {duality principle for biring}%
\Index{texBiring}
   {duality principle for biring of matrices}%
   {duality principle for biring of matrices}%
\SetIndexSpace%
\Index{texTstarMorphism}
   {effective representation of algebra $A$}%
   {effective representation of algebra}%
\Index{texFiberedAlgebra}
   {effective representation of fibered $\mathcal{F}$\Hyph algebra}%
   {effective representation of fibered F-algebra}%
\Index{texFiberedGroup}
   {effective \Ts representation of fibered group}%
   {effective representation of fibered group}%
\Index{texRepresentation}
   {effective representation of group}%
   {effective representation of group}%
\Index{texVectorSpace}
   {effective representation of skew field}%
   {effective representation of skew field}%
\Index{texELie}
   {enhanced Lie group}%
   {enhanced Lie group}%
\Index{texDiffEq}
   {essential parameters}%
   {essential parameters}%
\Index{texBundleRelation}
   {extension of correspondence}%
   {extension of correspondence}%
\Index{texAffine}
   {extreme line}%
   {extreme line}%
\SetIndexSpace%
\Index{texBundleRelation}
   {fibered correspondence from $\mathcal{A}$ to $\mathcal{B}$}%
   {fibered correspondence from A to B}%
\Index{texBundleRelation}
   {fibered correspondence in $\mathcal{A}$}%
   {fibered correspondence in A}%
\Index{texBundleRelation}
   {fibered correspondence of homomorphism}%
   {fibered correspondence of homomorphism}%
\Index{texBundleRelation}
   {fibered equivalence}%
   {fibered equivalence}%
\Index{texFiberedAlgebra}
   {fibered $\mathcal{F}$\Hyph algebra}%
   {fibered F-algebra}%
\Index{texFiberedAlgebra}
   {fibered $\mathcal{F}$\Hyph subalgebra}%
   {fibered F-subalgebra}%
\Index{texFiberedAlgebra}
   {fibered group}%
   {fibered group}%
\Index{texFiberedMorphism}
   {fibered identification morphism}%
   {fibered identification morphism}%
\Index{texFiberedMorphism}
   {fibered little group}%
   {fibered little group}%
\Index{texFiberedMorphism}
   {fibered natural morphism}%
   {fibered natural morphism}%
\Index{texBundleRelation}
   {fibered ordering}%
   {fibered ordering}%
\Index{texBundleRelation}
   {fibered preordering}%
   {fibered preordering}%
\Index{texFiberedAlgebra}
   {fibered ring}%
   {fibered ring}%
\Index{texFiberedMorphism}
   {fibered stability group}%
   {fibered stability group}%
\Index{texBundleRelation}
   {fibered subset}%
   {fibered subset}%
\Index{texNewton}
   {field-strength tensor}%
   {field-strength tensor}%
\Index{texBundleRelation}
   {filter $\mathfrak{F}$ converges to $A$}%
   {filter converges}%
\Index{texNewton}
   {first Newton law}%
   {First Newton law}%
\Index{texFiberedMorphism}
   {free \Ts representation of fibered group}%
   {free representation of fibered group}%
\Index{texTstarRepresentation}
   {free \Ts representation of group}%
   {free representation of group}%
\Index{texAffine}
   {Frenet transport}%
   {Frenet transport}%
\Index{texCalculus}
   {function continuous with respect to set of arguments}%
   {function continuous with respect to set of arguments}%
\Index{texCalculus}
   {function of $\gi n$ $D$\Hyph valued variables}%
   {function of n D valued variables}%
\SetIndexSpace%
\Index{texTypeBasis}
   {\Gbasis}%
   {G-basis}%
\Index{texBasis}
   {\Gbasis\ of vector space}%
   {G-basis}%
\Index{texBasis}
   {\Gcoords\ of basis}%
   {G-coordinates}%
\Index{texTypeBasis}
   {\Gcoords}%
   {G-coordinates}%
\Index{texTypeBasis}
   {\Gspace}%
   {GSpace}%
\Index{texDrcBasis}
   {geometrical object}%
   {}%
\SubIndex{texDrcBasis}
   {defined in vector space}%
   {geometrical object, vector space}%
\SubIndex{texDrcBasis}
   {in coordinate representation defined in vector space}%
   {geometrical object, coordinate vector space}%
\SubIndex{texDrcBasis}
   {of type $A$}%
   {geometrical object of type A, vector space}%
\Index{texBasis}
   {geometrical object in coordinate representation}%
   {geometrical object, coordinate vector space}%
\Index{texBasis}
   {geometrical object in vector space}%
   {geometrical object, vector space}%
\Index{texBasis}
   {geometrical object of type $A$ in vector space}%
   {geometrical object of type A, vector space}%
\Index{texGroupRing}
   {group algebra}%
   {group algebra}%
\Index{texBasis}
   {\Gspace}%
   {GSpace}%
\SetIndexSpace%
\Index{texBiring}
   {Hadamard inverse of matrix}%
   {Hadamard inverse of matrix}%
\Index{texReferenceFrame}
   {holonomic coordinates of vector}%
   {vector holonomic coordinates}%
\Index{texFiberedGroup}
   {homogeneous bundle of fibered group}%
   {homogeneous bundle of fibered group}%
\Index{texTstarRepresentation}
   {homogeneous space of group}%
   {homogeneous space of group}%
\Index{texRepresentation}
   {homogeneous space of group}%
   {homogeneous space of group}%
\Index{texFiberedAlgebra}
   {homomorphism of fibered $\mathcal{F}$\Hyph algebras}%
   {homomorphism of fibered F-algebras}%
\Index{texFiberedGroup}
   {homomorphism of fibered groups}%
   {homomorphism of fibered groups}%
\SetIndexSpace%
\Index{texLieRepresentation}
   {infinitesimal generator}%
   {infinitesimal generator}%
\Index{texLinearLie}
   {infinitesimal generators of group Lie}%
   {infinitesimal generators of group Lie}%
\Index{texDrcBasis}
   {invariance principle}%
   {invariance principle}%
\Index{texBasis}
   {invariance principle in vector space}%
   {invariance principle, vector space}%
\Index{texBundleRelation}
   {inverse fibered correspondence}%
   {inverse fibered correspondence}%
\Index{texBundleRelation}
   {inverse reduced fibered correspondence}%
   {inverse reduced fibered correspondence}%
\Index{texFiberedAlgebra}
   {isomorphism of fibered $\mathcal{F}$\Hyph algebras}%
   {isomorphism of fibered F-algebras}%
\SetIndexSpace%
\Index{texFiberedGroup}
   {kernel of inefficiency of representation of fibered group}%
   {kernel of inefficiency of representation of fibered group}%
\Index{texRepresentation}
   {kernel of inefficiency of representation of group}%
   {kernel of inefficiency of representation of group}%
\Index{texTstarRepresentation}
   {kernel of inefficiency of \Ts representation of group $G$}%
   {kernel of inefficiency of representation of group}%
\Index{texDiffProperty}
   {Killing equation}%
   {Killing equation}%
\Index{texDiffProperty}
   {Killing equation of second type}%
   {Killing equation second type}%
\Index{texDiffProperty}
   {Killing vector of second type}%
   {Killing vector second type}%
\Index{texBiring}
   {Kronecker symbol}%
   {Kronecker symbol}%
\SetIndexSpace%
\Index{texLie}
   {left invariant vector field}%
   {left invariant vector}%
\Index{texVectorSpace}
   {left module}%
   {left module}%
\Index{texTstarRepresentation}
   {left shift}%
   {left shift}%
\Index{texFiberedGroup}
   {left shift on fibered group}%
   {Tstar shift, fibered group}%
\Index{texRepresentation}
   {left shift on group}%
   {left shift, group}%
\Index{texLie}
   {left structural constant of Lie algebra}%
   {left structural constant of Lie algebra}%
\Index{texVectorSpace}
   {left vector space}%
   {left vector space}%
\Index{texRepresentation}
   {left-side contravariant representation of group}%
   {left-side contravariant representation}%
\Index{texRepresentation}
   {left-side covariant representation of group}%
   {left-side covariant representation}%
\Index{texTstarMorphism}
   {left-side representation of $\mathcal{F}$\Hyph algebra $A$ in set $M$}%
   {left-side representation of algebra}%
\Index{texFiberedAlgebra}
   {left-side representation of fibered $\mathcal{F}$\Hyph algebra}%
   {left-side representation of fibered F-algebra}%
\Index{texRepresentation}
   {left-side representation of group}%
   {left-side representation of group}%
\Index{texTstarMorphism}
   {left-side transformation}%
   {left-side transformation}%
\Index{texFiberedAlgebra}
   {left-side transformation on bundle}%
   {left-side transformation of bundle}%
\Index{texLie}
   {Lie algebra of Lie group}%
   {algebra Lie group Lie}%
\SubIndex{texLie}
   {left defined}%
   {left defined algebra Lie}%
\SubIndex{texLie}
   {right defined}%
   {right defined algebra Lie}%
\Index{texDiffProperty}
   {Lie derivative}%
   {Lie derivative}%
\SubIndex{texDiffProperty}
   {of connection}%
   {Lie derivative of connection}%
\SubIndex{texDiffProperty}
   {of metric}%
   {Lie derivative of metric}%
\Index{texLie}
   {Lie group basic operators}%
   {Lie group basic operators}%
\Index{texBundleRelation}
   {lift of correspondence}%
   {lift of correspondence}%
\Index{texCartesian}
   {lift of map}%
   {lift of map}%
\Index{texBundleRelation}
   {limit of correspondence with respect to the filter}%
   {limit of correspondence with respect to the filter}%
\Index{texBundleRelation}
   {limit of filter}%
   {limit of filter}%
\Index{texBundleRelation}
   {limit set of filter}%
   {limit set of filter}%
\Index{texRepresentation}
   {linear representation of group}%
   {linear representation of group}%
\Index{texReferenceFrame}
   {linearly dependent vector fields}%
   {linearly dependent vector fields}%
\Index{texTstarRepresentation}
   {little group}%
   {little group}%
\Index{texReferenceFrame}
   {local reference frame}%
   {local reference frame}%
\Index{texGenRelativity}
   {Lorentz transformation}%
   {Lorentz transformation}%
\SetIndexSpace%
\Index{texReferenceFrame}
   {map of type $G$ on manifold}%
   {map of type G on manifold}%
\Index{texCartesian}
   {mapping space}%
   {mapping space}%
\Index{texDrcMorphism}
   {matrix of \drc linear map}%
   {matrix of drc linear map}%
\Index{texAffine}
   {metric-affine manifold}%
   {metric-affine manifold}%
\Index{texTstarMorphism}
   {morphism of \Ts representations from $f$ into $g$}%
   {morphism of representations from f into g}%
\Index{texTstarMorphism}
   {morphism of \Ts representations of $\mathcal{F}$ algebra}%
   {morphism of representations of F algebra}%
\Index{texTstarMorphism}
   {morphism of \Ts representations of $\mathcal{F}$\Hyph algebra in $\mathcal{H}$\Hyph algebra}%
   {morphism of representations of F algebra in H algebra}%
\Index{}
   {movement on basis manifold}%
   {movement transformation}%
\SetIndexSpace%
\Index{texPolymodule}
   {$(n)$\hyph vector space}%
   {(n)-vector space}%
\Index{texBundleRelation}
   {$n$\Hyph ary fibered relation}%
   {fibered relation}%
\Index{texAffine}
   {nonmetricity}%
   {nonmetricity}%
\Index{texVectorSpace}
   {nonsingular \drc system of linear equations}%
   {nonsingular system of linear equations}%
\Index{texRepresentation}
   {nonsingular \Ts transformation}%
   {nonsingular transformation}%
\Index{texCalculus}
   {norm on \drc vector space}%
   {norm on drc vector space}%
\Index{texCalculus}
   {normed \drc vector space}%
   {normed drc vector space}%
\SetIndexSpace%
\Index{texFiberedAlgebra}
   {operation on bundle}%
   {operation on bundle}%
\Index{texBundleRelation}
   {opposite fibered preordering}%
   {opposite fibered preordering}%
\Index{texFiberedGroup}
   {orbit of representation of fibered group}%
   {orbit of representation of fibered group}%
\Index{texRepresentation}
   {orbit of representation of group}%
   {orbit of representation of group}%
\Index{texTstarRepresentation}
   {orbit of \Ts representation of group}%
   {orbit of representation of group}%
\Index{texBasis}
   {orthonornal basis}%
   {Orthonornal Basis}%
\SetIndexSpace%
\Index{texReferenceFrame}
   {parallelogram}%
   {parallelogram}%
\Index{texCalculus}
   {partial derivative of mapping $f$ with respect to variable $v^{\gi a}$}%
   {partial derivative of mapping with respect to variable, skew field}%
\Index{texCalculus}
   {partial derivative of mapping $\Vector f$ with respect to variable $v^{\gi a}$}%
   {partial derivative of mapping with respect to variable, drc vector space}%
\Index{texBasis}
   {passive representation}%
   {passive representation}%
\Index{texBasis}
   {passive transformation on basis manifold}%
   {passive transformation}%
\Index{texDrcBasis}
   {passive \Ts representation}%
   {passive representation}%
\Index{texReferenceFrame}
   {pfaffian derivative}%
   {pfaffian derivative}%
\Index{texNewton}
   {potential energy}%
   {potential energy}%
\Index{texDrcBasis}
   {product of geometrical object and constant}%
   {product of geometrical object and constant}%
\Index{texBasis}
   {product of geometrical object and constant in vector space}%
   {product of geometrical object and constant, vector space}%
\Index{texTstarMorphism}
   {product of morphisms of \Ts representations of $\mathcal{F}$\Hyph algebra}%
   {product of morphisms of representations of F algebra}%
\SetIndexSpace%
\Index{texBasis}
   {quasi affine transformation on basis manifold}%
   {quasi affine transformation}%
\Index{texBasis}
   {quasi movement on basis manifold}%
   {quasi movement}%
\Index{texFiberedMorphism}
   {quotient bundle}%
   {quotient bundle}%
\SetIndexSpace%
\Index{texBiring}
   {\sups row of matrix}%
   {r row}%
\Index{texVectorSpace}
   {\sups rows \drc vector space}%
   {sups rows drc vector space}%
\Index{texBiring}
   {$r$\hyph row of matrix}%
   {r-row}%
\Index{texBiring}
   {\RC inverse element of biring}%
   {rc-inverse element}%
\Index{texVectorSpace}
   {\RC major minor}%
   {RC-major minor}%
\Index{texVectorSpace}
   {\RC matrix group}%
   {rc-matrix group}%
\Index{texVectorSpace}
   {\RC nonsingular matrix}%
   {RC nonsingular matrix}%
\Index{texBiring}
   {\RC power}%
   {rc power}%
\Index{texBiring}
   {\RC product of matrices}%
   {rc-product of matrices}%
\Index{texBiring}
   {\RC quasideterminant}%
   {RC-quasideterminant}%
\Index{texVectorSpace}
   {\RC rank of matrix}%
   {rc-rank of matrix}%
\Index{texVectorSpace}
   {\RC singular matrix}%
   {RC singular matrix}%
\Index{texVectorSpace}
   {\rcd vector space}%
   {rcd vector space}%
\Index{texCartesian}
   {reduced Cartesian product of bundles}%
   {reduced Cartesian product of bundles}%
\Index{texCartesian}
   {reduced Cartesian product of total spaces}%
   {reduced Cartesian product of total spaces}%
\Index{texBundleRelation,texBundleRelation}
   {reduced fibered correspondence from $\mathcal{A}$ to $\mathcal{B}$}%
   {reduced fibered correspondence from A to B}%
\Index{texBundleRelation}
   {reduced fibered correspondence in $\mathcal{A}$}%
   {reduced fibered correspondence in A}%
\Index{texBiring}
   {reducible biring}%
   {reducible biring}%
\Index{texReferenceFrame}
   {reference frame}%
   {reference frame}%
\Index{texGenRelativity}
   {reference frame in event space}%
   {reference frame in event space}%
\Index{texReferenceFrame}
   {reference frame manifold}%
   {reference frame manifold}%
\Index{texBundleRelation}
   {reflexive $2$\Hyph ary fibered relation}%
   {reflexive 2 ary fibered relation}%
\Index{texTstarRepresentation}
   {representation of group}%
   {}%
\SubIndex{texTstarRepresentation}
   {contravariant \Ts}%
   {Tstar contravariant representation of group}%
\SubIndex{texTstarRepresentation}
   {covariant \Ts}%
   {Tstar covariant representation of group}%
\SubIndex{texDrcBasis}
   {\drc linear \sT}%
   {linear representation of group}%
\SubIndex{texTstarRepresentation}
   {effective}%
   {effective representation of group}%
\SubIndex{texDrcBasis}
   {\rcd}%
   {rcd linear representation of group}%
\SubIndex{texTstarRepresentation}
   {\sT}%
   {starT representation of group}%
\SubIndex{texTstarRepresentation}
   {\Ts}%
   {Tstar representation of group}%
\Index{texRepresentation}
   {representation of group}%
   {representation of group}%
\Index{texDrcBasis}
   {representative of geometrical object in vector space}%
   {representative of geometrical object, vector space}%
\Index{texBasis}
   {representative of geometrical object in vector space}%
   {representative of geometrical object, vector space}%
\Index{texBundleRelation}
   {restriction of correspondence $\Phi$ to set $C$}%
   {restriction of correspondence}%
\Index{texLie}
   {right invariant vector field}%
   {right invariant vector}%
\Index{texVectorSpace}
   {right module}%
   {right module}%
\Index{texTstarRepresentation}
   {right shift}%
   {right shift}%
\Index{texRepresentation}
   {right shift on group}%
   {right shift, group}%
\Index{texLie}
   {right structural constant of Lie algebra}%
   {right structural constant of Lie algebra}%
\Index{texVectorSpace}
   {right vector space}%
   {right vector space}%
\Index{texRepresentation}
   {right-side contravariant representation of group}%
   {right-side contravariant representation}%
\Index{texRepresentation}
   {right-side covariant representation of group}%
   {right-side covariant representation}%
\Index{texTstarMorphism}
   {right-side representation of $\mathcal{F}$\Hyph algebra $A$ in set $M$}%
   {right-side representation of algebra}%
\Index{texFiberedAlgebra}
   {right-side representation of fibered $\mathcal{F}$\Hyph algebra}%
   {right-side representation of fibered F-algebra}%
\Index{texRepresentation}
   {right-side representation of group}%
   {right-side representation of group}%
\Index{texTstarMorphism}
   {right-side transformation}%
   {right-side transformation}%
\Index{texRepresentation}
   {right-side transformation}%
   {right-side transformation}%
\Index{texRepresentation}
   {row vector}%
   {row vector}%
\Index{texVectorSpace}
   {$R\star$\Hyph module}%
   {Rstar-module}%
\SetIndexSpace%
\Index{texNewton}
   {scalar potential}%
   {scalar potential}%
\Index{texNewton}
   {second Newton law}%
   {Second Newton law}%
\Index{texTstarMorphism}
   {single transitive representation of algebra $A$}%
   {single transitive representation of algebra}%
\Index{texFiberedAlgebra}
   {single transitive representation of fibered $\mathcal{F}$\Hyph algebra}%
   {single transitive representation of fibered F-algebra}%
\Index{texRepresentation}
   {single transitive representation of group}%
   {single transitive representation of group}%
\Index{texTstarRepresentation}
   {space of orbits of \Ts representation}%
   {space of orbits of Ts representation}%
\Index{texTidal}
   {speed of deviation}%
   {speed of deviation}%
\Index{texDrcMorphism}
   {$(S\RCstar,T\RCstar)$\Hyph linear map of vector spaces}%
   {src trc linear map of vector spaces}%
\Index{texTstarRepresentation}
   {stability group}%
   {stability group}%
\Index{texDrcBasis}
   {standard coordinates of basis}%
   {standard coordinates of basis}%
\Index{texBasis}
   {standard coordinates of basis}%
   {standard coordinates of basis}%
\Index{texBiring}
   {standard representation of matrix}%
   {Standard representation}%
\Index{texVectorSpace}
   {$\star D$\hyph vector space}%
   {starD-vector space}%
\Index{texLinearMap}
   {$\star D$\Hyph product of \drc linear map $A$ over scalar}%
   {starD product of drc linear map over scalar}%
\Index{texVectorSpace}
   {$\star R$\hyph module}%
   {starR-module}%
\Index{texTstarRepresentation}
   {\sT shift}%
   {starT shift}%
\Index{texFiberedGroup}
   {\sT shift on fibered group}%
   {starT shift, fibered group}%
\Index{texTstarMorphism}
   {\sT representation of $\mathcal{F}$\Hyph algebra $A$ in set $M$}%
   {starT representation of algebra}%
\Index{texFiberedAlgebra}
   {\sT representation of fibered $\mathcal{F}$\Hyph algebra}%
   {starT representation of fibered F-algebra}%
\Index{texFiberedGroup}
   {\sT representation of fibered group}%
   {starT representation of fibered group}%
\Index{texTstarMorphism}
   {\sT transformation}%
   {starT transformation}%
\Index{texFiberedAlgebra}
   {\sT transformation on bundle}%
   {starT transformation of bundle}%
\SubIndex{}
   {nonsingular}%
   {nonsingular transformation of bundle}%
\Index{texBundleRelation}
   {subbundle}%
   {subbundle}%
\Index{texLinearMap}
   {sum of \drc linear maps}%
   {sum of drc linear maps, drc vector spaces}%
\Index{texDrcBasis}
   {sum of geometrical objects}%
   {sum of geometrical objects}%
\Index{texBasis}
   {sum of geometrical objects in vector space}%
   {sum of geometrical objects, vector space}%
\Index{texBundleRelation}
   {symmetric $2$\Hyph ary fibered relation}%
   {symmetric 2 ary fibered relation}%
\Index{texBasis}
   {symmetry group}%
   {symmetry group}%
\Index{texDrcBasis}
   {symmetry group}%
   {SymmetryGroup}%
\Index{texGenRelativity}
   {synchronization of reference frame}%
   {synchronization of reference frame}%
\SetIndexSpace%
\Index{texLie}
   {tensor product of representations}%
   {tensor product of representations}%
\Index{texCalculus}
   {topological \drc vector space}%
   {topological drc vector space}%
\Index{texCalculus}
   {topological skew field}%
   {topological skew field}%
\Index{texAffine}
   {torsion}%
   {torsion}%
\Index{texFiberedMorphism}
   {tower of bundles}%
   {tower of bundles}%
\Index{texTstarMorphism}
   {transformation of set}%
   {transformation of set}%
\Index{texDrcBasis}
   {transformation on basis manifold}%
   {}%
\SubIndex{texDrcBasis}
   {active}%
   {active transformation, vector space}%
\SubIndex{texTypeBasis}
   {affine}%
   {affine transformation}%
\SubIndex{texTypeBasis}
   {movement}%
   {movement transformation}%
\SubIndex{texDrcBasis}
   {passive}%
   {passive transformation, vector space}%
\SubIndex{texTypeBasis}
   {quasi affine}%
   {quasi affine transformation}%
\SubIndex{texTypeBasis}
   {quasi movement}%
   {quasi movement}%
\Index{texFiberedAlgebra}
   {transformation on bundle}%
   {transformation of bundle}%
\Index{texBundleRelation}
   {transitive $2$\Hyph ary fibered relation}%
   {transitive 2 ary fibered relation}%
\Index{texTstarMorphism}
   {transitive representation of algebra $A$}%
   {transitive representation of algebra}%
\Index{texFiberedAlgebra}
   {transitive representation of fibered $\mathcal{F}$\Hyph algebra}%
   {transitive representation of fibered F-algebra}%
\Index{texRepresentation}
   {transitive representation of group}%
   {transitive representation of group}%
\Index{texVectorSpace}
   {\Ts linear composition of  vectors}%
   {linear composition of  vectors}%
\Index{texVectorSpace}
   {\Ts matrices vector space}%
   {matrices vector space}%
\Index{texTstarRepresentation}
   {\Ts shift}%
   {Tstar shift}%
\Index{texTstarMorphism}
   {\Ts representation of $\mathcal{F}$\Hyph algebra $A$ in set $M$}%
   {Tstar representation of algebra}%
\Index{texFiberedAlgebra}
   {\Ts representation of fibered $\mathcal{F}$\Hyph algebra}%
   {Tstar representation of fibered F-algebra}%
\Index{texFiberedGroup}
   {\Ts representation of fibered group}%
   {Tstar representation of fibered group}%
\Index{texTstarMorphism}
   {\Ts transformation}%
   {Tstar transformation}%
\Index{texFiberedAlgebra}
   {\Ts transformation on bundle}%
   {Tstar transformation of bundle}%
\Index{texFiberedGroup}
   {twin representations of fibered group}%
   {twin representations of fibered group}%
\Index{texTstarRepresentation}
   {twin representations of group}%
   {twin representations of group}%
\Index{texLinearMap}
   {twin representations of skew field}%
   {twin representations of skew field}%
\Index{texReferenceFrame}
   {type $G$ reference frame}%
   {type G reference frame}%
\SetIndexSpace%
\Index{texVectorSpace}
   {unitarity law}%
   {}%
\SubIndex{texVectorSpace}
   {for \Ts vector space}%
   {unitarity law, Tstar vector space}%
\SetIndexSpace%
\Index{texPolymodule}
   {($D_1\RCstar$, ..., $D_n\RCstar$)\hyph vector space}%
   {(d1rc,dnrc)-vector space}%
\Index{texPolymodule}
   {($S\star$, $\star T$)\hyph vector space}%
   {(Sstar,starT)-vector space}%
\Index{texCalculus}
   {valued skew field}%
   {valued skew field}%
\Index{texFiberedAlgebra}
   {vector bundle}%
   {vector bundle}%
\Index{texNewton}
   {vector potential}%
   {vector potential}%
\Index{texVectorSpace}
   {vector space type}%
   {vector space type}%

\CloseIndex

\def\indexname{Special Symbols and Notations}
\OpenIndex

\SetIndexSpace%
\Symb{texBiring}
   {$(^{\gi a}_{\gi b})$\hyph\CR quasideterminant}%
   {a b CR quasideterminant definition}%
\Symb{texBiring}
   {$(^{\gi a}_{\gi b})$\hyph \RC quasideterminant}%
   {a b RC-quasideterminant definition}%
\Symb{texBiring}
   {minor}%
   {A from b a}%
\Symb{texBiring}
   {minor}%
   {A from columns T}%
\Symb{texBiring}
   {minor}%
   {A from rows S}%
\Symb{texBiring}
   {minor}%
   {A without column a}%
\Symb{texBiring}
   {minor}%
   {A without columns T}%
\Symb{texBiring}
   {minor}%
   {A without row b}%
\Symb{texBiring}
   {minor}%
   {A without rows S}%
\Symb{texPolymodule}
   {active transformation}%
   {active transformation}%
\Symb{texTypeBasis}
   {affine space}%
   {affine space}%
\Symb{texBasis}
   {affine space}%
   {An}%
\Symb{texBiring}
   {\subs row ($c$\hyph row) of matrix}%
   {c row}%
\Symb{texBiring}
   {\CR power of element $A$ of biring}%
   {cr power}%
\Symb{texBiring}
   {\CR inverse element of biring}%
   {cr-inverse element}%
\Symb{texBiring}
   {\CR product of matrices}%
   {cr-product of matrices}%
\Symb{texVectorSpace}
   {\dcr vector}%
   {dcr vector}%
\Symb{texLie}
   {derivative of left shift}%
   {derivative of left shift}%
\Symb{texLie}
   {derivative of left shift}%
   {derivative of left shift, 1-Parameter Group}%
\Symb{texLie}
   {derivative of right shift}%
   {derivative of right shift}%
\Symb{texLie}
   {}%
   {derivative of right shift}%
\Symb{texLie}
   {derivative of right shift}%
   {derivative of right shift, 1-Parameter Group}%
\Symb{texLie}
   {derivative of left shift}%
   {derivative of Tstar shift}%
\Symb{texVectorSpace}
   {\drc vector}%
   {drc vector}%
\Symb{texAffine}
   {derivative}%
   {overline nabla_l, definition 2}%
\Symb{texPolymodule}
   {passive transformation}%
   {passive transformation}%
\Symb{texBiring}
   {\sups row ($r$\hyph row) of matrix}%
   {r row}%
\Symb{texBiring}
   {\RC power of element $A$ of biring}%
   {rc power}%
\Symb{texBiring}
   {\RC inverse element of biring}%
   {rc-inverse element}%
\Symb{texBiring}
   {\RC product of matrices}%
   {rc-product of matrices}%
\Symb{texBiring}
   {\RC quasideterminant}%
   {RC-quasideterminant definition}%
\Symb{texTstarRepresentation}
   {right shift}%
   {right shift}%
\Symb{texFiberedGroup}
   {\sT shift}%
   {starT shift, fibered group}%
\Symb{texTstarRepresentation}
   {left shift}%
   {Tstar shift}%
\Symb{texFiberedGroup}
   {\Ts shift}%
   {Tstar shift, fibered group}%
\Symb{texReferenceFrame}
   {anholonomic coordinates of vector}%
   {vector anholonomic coordinates}%
\Symb{texReferenceFrame}
   {holonomic coordinates of vector}%
   {vector holonomic coordinates}%

\SetIndexSpace%
\Symb{texBasis}
   {basis manifold of affine space}%
   {BAn}%
\Symb{texReferenceFrame}
   {basis manifold of manifold}%
   {basis manifold of manifold}%
\Symb{texPolymodule}
   {basis manifold of vector space}%
   {basis manifold of vector space}%
\Symb{texBasis}
   {basis manifold of vector space $\mathcal{V}$}%
   {basis manifold of vector space}%
\Symb{texBasis}
   {basis manifold of central affine space}%
   {BCAn}%
\Symb{texBasis}
   {basis manifold of Euclid space}%
   {BEn}%
\Symb{texCartesian}
   {Cartesian power $A$ of set $B$}%
   {Cartesian power of set}%
\Symb{texTypeBasis}
   {basis manifold of affine space}%
   {FAn}%
\Symb{texTypeBasis}
   {basis manifold of central affine space}%
   {FCAn}%
\Symb{texTypeBasis}
   {basis manifold of Euclid space}%
   {FEn}%

\SetIndexSpace%
\Symb{texBasis}
   {central affine space}%
   {CAn}%
\Symb{texTypeBasis}
   {central affine space}%
   {central affine space}%
\Symb{texLie}
   {left structural constant of Lie algebra}%
   {left structural constant of Lie algebra}%
\Symb{texLie}
   {right structural constant of Lie algebra}%
   {right structural constant of Lie algebra}%
\Symb{texReferenceFrame}
   {set of functions defined on manifold}%
   {set of functions defined on manifold}%

\SetIndexSpace%
\Symb{texLieRepresentation}
   {basis vector of \sT representation}%
   {basis vector of starT representation}%
\Symb{texLieRepresentation}
   {basis vector of \sT representation}%
   {basis vector of starT representation, coordinates}%
\Symb{texLieRepresentation}
   {basis vector of \Ts representation}%
   {basis vector of Tstar representation}%
\Symb{texLieRepresentation}
   {basis vector of \Ts representation}%
   {basis vector of Tstar representation, coordinates}%
\Symb{texVectorSpace}
   {\subs rows \dcr vector space}%
   {c rows dcr vector space}%
\Symb{texCalculus}
   {differential of function}%
   {differential, drc vector space to drc vector space}%
\Symb{texCalculus}
   {differential of function}%
   {differential, drc vector space to skew field}%
\Symb{texVectorSpace}
   {\drc coordinate vector space}%
   {drc coordinate vector space}%
\Symb{texVectorSpace}
   {matrices vector space}%
   {matrices vector space}%
\Symb{texAffine}
   {Cartan derivative}%
   {overbrace D}%
\Symb{texAffine}
   {derivative}%
   {overline D}%
\Symb{texCalculus}
   {partial derivative of mapping $\Vector f$ with respect to variable $v^{\gi a}$}%
   {partial derivative of mapping, 1, drc vector space}%
\Symb{texCalculus}
   {partial derivative of mapping $f$ with respect to variable $v^{\gi a}$}%
   {partial derivative of mapping, 1, skew field}%
\Symb{texVectorSpace}
   {\sups rows \drc vector space}%
   {r rows drc vector space}%
\Symb{texTidal}
   {speed of deviation}%
   {speed of deviation}%
\Symb{texVectorSpace}
   {vector space type}%
   {vector space type}%

\SetIndexSpace%
\Symb{texTypeBasis}
   {affine basis}%
   {Affine Basis}%
\Symb{texBasis}
   {affine basis}%
   {Affine Basis}%
\Symb{texTypeBasis}
   {basis}%
   {basis}%
\Symb{texBasis}
   {basis of vector space}%
   {Basis e}%
\Symb{texBasis}
   {basis in vector space $\mathcal{V}$}%
   {basis in V}%
\Symb{texVectorSpace}
   {basis of vector space}%
   {basis, vector space}%
\Symb{texPolymodule}
   {basis of $(n)$\hyph vector space}%
   {basis,n vector space}%
\Symb{texCartesian}
   {Cartesian power of total spaces}%
   {Cartesian power of total spaces}%
\Symb{texCartesian}
   {Cartesian product of total spaces}%
   {Cartesian product of total spaces, definition 1}%
\Symb{texBasis}
   {central affine basis}%
   {Central Affine Basis}%
\Symb{texReferenceFrame}
   {form of reference frame}%
   {dual forms, reference frame}%
\Symb{texTypeBasis}
   {Euclid space}%
   {En}%
\Symb{texBasis}
   {Euclid space}%
   {En}%
\Symb{texBasis}
   {pseudo Euclid space}%
   {Enm}%
\Symb{texTypeBasis}
   {pseudo Euclid space}%
   {Enm}%
\Symb{texFiberedAlgebra}
   {identical transformation of bundle}%
   {identical transformation of bundle}%
\Symb{texBasis}
   {orthonornal basis}%
   {Orthonornal Basis}%
\Symb{texCartesian}
   {reduced Cartesian product of total spaces}%
   {reduced Cartesian product of total spaces, definition 1}%
\Symb{texFiberedAlgebra}
   {set of nonsingular \sT transformations of bundle $\mathcal{E}$}%
   {set of starT nonsingular transformations of bundle}%
\Symb{texFiberedAlgebra}
   {set of nonsingular \Ts transformations of bundle $\mathcal{E}$}%
   {set of Tstar nonsingular transformations of bundle}%
\Symb{texBasis}
   {standard coordinates of basis}%
   {standard coordinates of basis}%
\Symb{texReferenceFrame}
   {standard coordinates of reference frame}%
   {standard coordinates of reference frame}%
\Symb{texReferenceFrame}
   {vector field of reference frame}%
   {vector field of reference frame}%
\Symb{texReferenceFrame}
   {vector field of reference frame}%
   {vector field, reference frame}%
\Symb{texBasis}
   {vector of basis}%
   {vector of basis}%

\SetIndexSpace%
\Symb{texVectorSpace}
   {coordinates of basis in \subs rows \dcr vector space}%
   {basis coordinates, c rows dcr vector space}%
\Symb{texVectorSpace}
   {coordinates of basis in \sups rows \drc vector space}%
   {basis coordinates, r rows drc vector space}%
\Symb{texVectorSpace}
   {basis for \subs rows \dcr vector space}%
   {basis, c rows dcr vector space}%
\Symb{texVectorSpace}
   {basis for \sups rows \drc vector space}%
   {basis, r rows drc vector space}%
\Symb{texDiffEq}
   {central affine basis}%
   {Central Affine Basis}%
\Symb{texReferenceFrame}
   {coordinate reference frame}%
   {coordinate reference frame}%
\Symb{texBundleRelation}
   {filter $\mathfrak{F}$ converges to set $A$}%
   {filter converges}%
\Symb{texFiberedAlgebra}
   {homomorphism of fibered $\mathcal{F}$\Hyph algebras}%
   {homomorphism of fibered F-algebras}%
\Symb{texBundleRelation}
   {inverse fibered correspondence}%
   {inverse fibered correspondence, 1}%
\Symb{texBundleRelation}
   {inverse reduced fibered correspondence}%
   {inverse reduced fibered correspondence, 1}%
\Symb{texCartesian}
   {map to Cartesian product}%
   {map to Cartesian product}%
\Symb{texTstarRepresentation}
   {representation orbit of group $G$}%
   {orbit of representation of group}%
\Symb{texTypeBasis}
   {orthonornal basis}%
   {Orthonornal Basis}%
\Symb{texReferenceFrame}
   {reference frame}%
   {reference frame}%
\Symb{texReferenceFrame}
   {reference frame}%
   {reference frame, extensive definition}%
\Symb{texPolymodule}
   {standard coordinates of basis}%
   {standard coordinates of basis}%
\Symb{texPolymodule}
   {vector of basis}%
   {vector of basis}%

\SetIndexSpace%
\Symb{texVectorSpace}
   {\CR matrix group}%
   {cr-matrix group}%
\Symb{texFiberedMorphism}
   {fibered little group of section $h$}%
   {fibered little group}%
\Symb{texFiberedMorphism}
   {fibered stability group of section $h$}%
   {fibered stability group}%
\Symb{texLie}
   {algebra Lie of group Lie}%
   {g}%
\Symb{texLie}
   {left defined algebra Lie of group Lie}%
   {gl}%
\Symb{texTypeBasis}
   {affine transformation group}%
   {GLAn}%
\Symb{texBasis}
   {affine transformation group}%
   {GLAn}%
\Symb{texLie}
   {right defined algebra Lie of group Lie}%
   {gr}%
\Symb{texBasis}
   {group of homomorphisms of vector space $\mathcal{V}$}%
   {GV}%
\Symb{texTstarRepresentation}
   {little group of $x$}%
   {little group}%
\Symb{texFiberedGroup}
   {orbit of effective covariant \sT representation of fibered group}%
   {orbit of effective starT covariant representation of fibered group}%
\Symb{texTstarRepresentation}
   {orbit of effective covariant \sT representation of group}%
   {orbit of effective starT covariant representation of group}%
\Symb{texFiberedGroup}
   {orbit of effective covariant \Ts representation of fibered group}%
   {orbit of effective Tstar covariant representation of fibered group}%
\Symb{texTstarRepresentation}
   {orbit of effective covariant \Ts representation of group}%
   {orbit of effective Tstar covariant representation of group}%
\Symb{texVectorSpace}
   {\RC matrix group}%
   {rc-matrix group}%
\Symb{texTstarRepresentation}
   {stability group of $x$}%
   {stability group}%

\SetIndexSpace%
\Symb{texBiring}
   {Hadamard inverse of matrix}%
   {Hadamard inverse of matrix}%
\Symb{texLinearMap}
   {\rcd vector space of \drc linear maps}%
   {rcd vector space of drc linear maps}%

\SetIndexSpace%
\Symb{texLieRepresentation}
   {infinitesimal generator of representation}%
   {infinitesimal generator of representation}%
\Symb{texLinearLie}
   {Lie group infinitesimal generators}%
   {Lie group infinitesimal generators}%

\SetIndexSpace%
\Symb{texRepresentation}
   {left shift}%
   {left shift}%
\Symb{texDiffProperty}
   {Lie derivative of connection}%
   {Lie derivative of connection}%
\Symb{texDiffProperty}
   {Lie derivative of metric}%
   {Lie derivative of metric}%
\Symb{texBundleRelation}
   {limit of correspondence $\Phi$ with respect to the filter $\mathfrak{F}$}%
   {limit of correspondence with respect to the filter}%
\Symb{texBasis}
   {passive transformation}%
   {passive transformation}%
\Symb{texRepresentation}
   {set of left-side nonsingular transformations of set $M$}%
   {set of left-side nonsingular transformations}%

\SetIndexSpace%
\Symb{texTstarMorphism}
   {set of \sT transformations of set $M$}%
   {set of starT transformations}%
\Symb{texTstarMorphism}
   {set of \Ts transformations of set $M$}%
   {set of Tstar transformations}%
\Symb{texTstarRepresentation}
   {space of orbits of effective \sT covariant representation of the group}%
   {space of orbits of effective sT representation}%
\Symb{texTstarRepresentation}
   {space of orbits of effective \Ts covariant representation of the group}%
   {space of orbits of effective Ts representation}%
\Symb{texTstarRepresentation}
   {space of orbits of \Ts representation $f$ of group $G$ in set $M$}%
   {space of orbits of Ts representation}%

\SetIndexSpace%
\Symb{texBasis}
   {geometrical object in coordinate representation}%
   {geometrical object, coordinate vector space}%
\Symb{texBasis}
   {geometrical object}%
   {geometrical object, vector space}%
\Symb{texFiberedGroup}
   {orbit of representation of fibered group $\mathcal{G}$}%
   {orbit of representation of fibered group}%
\Symb{texRepresentation}
   {orbit of representation of the group $G$}%
   {orbit of representation of group}%

\SetIndexSpace%
\Symb{texCartesian}
   {bundle}%
   {bundle}%
\Symb{texFiberedMorphism}
   {bundle of level $2$}%
   {bundle of level 2}%
\Symb{texFiberedMorphism}
   {bundle of level $n$}%
   {bundle of level n}%
\Symb{texCartesian}
   {Cartesian power of bundle}%
   {Cartesian power of bundle}%
\Symb{texCartesian}
   {Cartesian product of bundles}%
   {Cartesian product of bundles, definition 1}%
\Symb{texCartesian}
   {reduced Cartesian product of bundles}%
   {reduced Cartesian product of bundles, definition 1}%
\Symb{texFiberedAlgebra}
   {set of nonsingular \sT transformations of bundle $\bundle{}pE{}$}%
   {set of starT nonsingular transformations of bundle, projection}%
\Symb{texFiberedAlgebra}
   {set of nonsingular \Ts transformations of bundle $\bundle{}pE{}$}%
   {set of Tstar nonsingular transformations of bundle, projection}%

\SetIndexSpace%
\Symb{texBasis}
   {active transformation}%
   {active transformation}%
\Symb{texAffine}
   {Cartan curvature}%
   {Cartan curvature}%
\Symb{texVectorSpace}
   {\CR rank of matrix}%
   {cr-rank of matrix}%
\Symb{texBundleRelation}
   {diagonal in bundle  $\bundle{}pA{}$}%
   {diagonal in bundle, 2}%
\Symb{texBundleRelation}
   {diagonal in bundle $\mathcal{A}$}%
   {diagonal in reduced bundle, 2}%
\Symb{texAffine}
   {curvature}%
   {GLn curvature_overline}%
\Symb{texVectorSpace}
   {\RC rank of matrix}%
   {rc-rank of matrix}%
\Symb{texRepresentation}
   {right shift}%
   {right shift}%
\Symb{texRepresentation}
   {set of right-side nonsingular transformations of set $M$}%
   {set of right-side nonsingular transformations}%

\SetIndexSpace%
\Symb{texBundleRelation}
   {composition of fibered correspondences}%
   {composition of fibered correspondences}%
\Symb{texBundleRelation}
   {inverse fibered correspondence}%
   {inverse fibered correspondence, 2}%
\Symb{texBundleRelation}
   {inverse reduced fibered correspondence}%
   {inverse reduced fibered correspondence, 2}%
\Symb{texVectorSpace}
   {linear span in vector space}%
   {linear span, vector space}%

\SetIndexSpace%
\Symb{texLie}
   {tangent plane to group $G$}%
   {TaG}%

\SetIndexSpace%
\Symb{texBasis}
   {coordinate vector space}%
   {coordinate vector space}%
\Symb{texBasis}
   {coordinates in vector space}%
   {coordinates in vector space}%
\Symb{texVectorSpace}
   {\dcr vector space}%
   {left CR vector space}%
\Symb{texVectorSpace}
   {\drc vector space}%
   {left RC vector space}%
\Symb{texLinearMap}
   {($S$, $T$)\hyph bimodule}%
   {R S bimodule}%
\Symb{texVectorSpace}
   {\crd vector space}%
   {right CR vector space}%
\Symb{texVectorSpace}
   {\rcd vector space}%
   {right RC vector space}%
\Symb{texBasis}
   {vector space}%
   {V}%
\Symb{texReferenceFrame}
   {vector space of vector fields}%
   {vector space of vector fields}%

\SetIndexSpace%
\Symb{texPolymodule}
   {geometrical object in coordinate representation		defined in vector space}%
   {geometrical object, coordinate vector space}%
\Symb{texPolymodule}
   {geometrical object in vector space}%
   {geometrical object, vector space}%

\SetIndexSpace%
\Symb{texReferenceFrame}
   {anholonomic coordinate}%
   {x(k)}%

\SetIndexSpace%
\Symb{texBundleRelation}
   {diagonal in bundle $\mathcal{A}$}%
   {diagonal in bundle, 1}%

\SetIndexSpace%
\Symb{texTidal}
   {deviation of trajectories}%
   {deviation of trajectories}%
\Symb{texRepresentation}
   {identical transformation}%
   {identical transformation}%
\Symb{texTstarMorphism}
   {identical transformation}%
   {identical transformation}%
\Symb{texBasis}
   {image of vector $\Vector e_k\in\Basis e$ under isomorphism to coordinate vector space}%
   {image of vector e_k, coordinate vector space}%
\Symb{texBiring}
   {Kronecker symbol}%
   {Kronecker symbol}%

\SetIndexSpace%
\Symb{texReferenceFrame}
   {anholonomic coordinates of connection}%
   {anholonomic coordinates of connection}%
\Symb{texAffine}
   {Cartan symbol}%
   {Cartan symbol}%
\Symb{texAffine}
   {connection}%
   {conection overline}%
\Symb{texAffine}
   {Cartan connection}%
   {overbrace Gamma i kl}%
\Symb{texCartesian}
   {set of sections of bundle}%
   {set of sections of bundle}%

\SetIndexSpace%
\Symb{texLie}
   {inverse operator to operator $\psi_l$}%
   {inverse operator to operator psi l}%
\Symb{texLie}
   {inverse operator to operator $\psi_r$}%
   {inverse operator to operator psi r}%

\SetIndexSpace%
\Symb{texReferenceFrame}
   {anholonomity object}%
   {anholonomity object}%

\SetIndexSpace%
\Symb{texLie}
   {basic operator of group Lie}%
   {Lie Basic Operator L}%
\Symb{texLie}
   {}%
   {Lie Basic Operator L}%
\Symb{texLie}
   {basic operator of group Lie}%
   {Lie Basic Operator L, 1-Parameter Group}%
\Symb{texLie}
   {basic operator of group Lie}%
   {Lie Basic Operator R}%
\Symb{texLie}
   {}%
   {Lie Basic Operator R}%
\Symb{texLie}
   {basic operator of group Lie}%
   {Lie Basic Operator R, 1-Parameter Group}%

\SetIndexSpace%
\Symb{texReferenceFrame}
   {coordinate reference frame}%
   {coordinate reference frame, extensive definition}%
\Symb{texCalculus}
   {partial derivative of mapping $\Vector f$ with respect to variable $v^{\gi a}$}%
   {partial derivative of mapping, 2, drc vector space}%
\Symb{texCalculus}
   {partial derivative of mapping $f$ with respect to variable $v^{\gi a}$}%
   {partial derivative of mapping, 2, skew field}%
\Symb{texReferenceFrame}
   {derivative $e_{(k)}$}%
   {partial(k)}%

\SetIndexSpace%
\Symb{texLie}
   {Lie group composition law}%
   {Lie group composition law}%

\SetIndexSpace%
\Symb{texAffine}
   {Cartan derivative}%
   {overbrace nabla_l}%
\Symb{texAffine}
   {}%
   {overbrace nabla_l}%
\Symb{texAffine}
   {derivative}%
   {overline nabla_l, definition 1}%

\SetIndexSpace%
\Symb{texBundleRelation}
   {restriction of correspondence $\Phi$ to set $C$}%
   {restriction of correspondence}%

\SetIndexSpace%
\Symb{texCartesian}
   {Cartesian product of bundles}%
   {Cartesian product of bundles, definition 2}%
\Symb{texCartesian}
   {Cartesian product of total spaces}%
   {Cartesian product of total spaces, definition 2}%
\Symb{texCartesian}
   {reduced Cartesian product of bundles}%
   {reduced Cartesian product of bundles, definition 2}%
\Symb{texCartesian}
   {reduced Cartesian product of total spaces}%
   {reduced Cartesian product of total spaces, definition 2}%

\SetIndexSpace%
\Symb{texBundleRelation}
   {fibered subset}%
   {fibered subset}%
\Symb{texBundleRelation}
   {subbundle}%
   {subbundle}%

\CloseIndex

\end{document}